\documentclass[reqno,10pt]{article}

\usepackage{amsmath,amsfonts,amsthm,mathrsfs,amssymb,bbm}%,eufrak}
\usepackage[all]{xy}
\usepackage{palatino}
\usepackage{a4wide}
\usepackage[%pdftex,
            breaklinks=true,
            pdftitle={Hecke Modifications, Wonderful Compactifications and Moduli of Principal Bundles},
            pdfauthor={Michael Lennox Wong}]{hyperref}
\usepackage{color}
\definecolor{tocolor}{rgb}{.1,.1,.5}
\definecolor{urlcolor}{rgb}{.2,.2,.6}
\definecolor{linkcolor}{rgb}{.1,.4,.6}
\definecolor{citecolor}{rgb}{.6,.3,.1}
\hypersetup{colorlinks=true, urlcolor=urlcolor, linkcolor=linkcolor, citecolor=citecolor}

\newcommand{\C}{\mathbb{C}}
\newcommand{\R}{\mathbb{R}}

\newcommand{\Z}{\mathbb{Z}}
\newcommand{\N}{\mathbb{N}}

\renewcommand{\1}{\mathbbm{1}}
\newcommand{\Res}{\text{Res}}

\newcommand{\supp}{\text{supp}}

\newcommand{\mero}{\text{mero}}

\newcommand{\arsim}{\xrightarrow{\sim}}

\newcommand{\G}{\mathscr{G}}
\newcommand{\Gr}{\text{Gr}}

\renewcommand{\b}{\mathfrak{b}}

\renewcommand{\L}{\mathscr{L}}

\newcommand{\NN}{\mathscr{N}}

\newcommand{\OO}{\mathscr{O}}
\renewcommand{\P}{\mathbb{P}}

\def\pd#1#2{\frac{\partial #1}{\partial #2}}

\newcommand{\g}{\mathfrak{g}}

\renewcommand{\t}{\mathfrak{t}}

\renewcommand{\sl}{\mathfrak{sl}}

\renewcommand{\u}{\mathfrak{u}}

\newcommand{\Ad}{\text{Ad}}
\newcommand{\ad}{\text{ad}}

\newcommand{\Aut}{\text{Aut}}
\newcommand{\End}{\text{End}}

\renewcommand{\L}{\mathscr{L}}
\renewcommand{\b}{\mathfrak{b}}
\newcommand{\GG}{\overline{G}}
\newcommand{\QQ}{\mathcal{Q}}

\renewcommand{\NN}{\mathscr{N}}

\newcommand{\Cz}{\C \{\!\{ z \}\!\}}
\newcommand{\Czz}{\C \{\!( z )\!\}}
\newcommand{\ev}{\text{ev}}
\newcommand{\af}{\text{af}}

\numberwithin{equation}{section}
\theoremstyle{definition}

\newtheorem{lem}[equation]{Lemma}
\newtheorem{thm}[equation]{Theorem}
\newtheorem{prop}[equation]{Proposition}
\newtheorem{cor}[equation]{Corollary}

\newtheorem{rmk}[equation]{Remark}

\begin{document}

\title{Hecke Modifications, Wonderful Compactifications and Moduli of Principal Bundles}

\author{Michael Lennox Wong\thanks{\tt wong@math.tifr.res.in}}
\date{}
%\address{School of Mathematics\\ Tata Institute of Fundamental Research\\ Dr.\ Homi Bhabha Road\\ Mumbai 400 005, India}
%\email{{\tt }}

% \subjclass[2010]{14D20 (Primary), 32G08 (Secondary)}

\maketitle

\begin{abstract}
In this paper, we obtain parametrizations of the moduli space of principal bundles over a compact Riemann surface using spaces of Hecke modifications in several cases.  We begin with a discussion of Hecke modifications for principal bundles and give constructions of ``universal'' Hecke modifications of a fixed bundle of fixed type.  This is followed by an overview of the construction of the ``wonderful,'' or De Concini--Procesi, compactification of a semi-simple algebraic group of adjoint type.  The compactification plays an important role in the deformation theory used in constructing the parametrizations.  A general outline to construct parametrizations is given and verifications for specific structure groups are carried out.
\end{abstract}

\section*{Introduction}

The main goal of this paper is to parametrize the moduli space of principal bundles over a compact Riemann surface using appropriate (symmetric) products of spaces of Hecke modifications of a fixed bundle.  A Hecke modification of a fixed bundle is obtained by ``twisting'' the transition function of that bundle near a point.  While neither the idea nor the application to moduli questions is new, the theory in the principal bundle setting does not seem to be well-developed and another goal here is to begin to fill in this lacuna.  The notion of a Hecke modification has its roots in Weil's concept of a ``matrix divisor'' \cite{Weil}, and is the basis of A.N.\ Tjurin's parametrization of the moduli space of rank $n$, degree $ng$ vector bundles over a Riemann surface of genus $g$ \cite{Tjurin}.  The notion has even further reach, as it is related to that of a Hecke operator acting on spaces of cusp forms (see \cite{HarderKazhdan}); thus, they play an important role in the geometrization of the Langlands program.

Tjurin's construction was later generalized to bundles of arbitrary degree by J.\ Hurtubise \cite{Hurtubise}.  The latter work makes use of the fact that $GL_n \C$ is open and dense in the space of $n \times n$ matrices, so for a more general structure group, one would like to embed the group as an open dense set in some larger space.  This is the entry point of the wonderful compactification.  Originally conceived to attack problems in enumerative geometry, this construction was first obtained by C.\ De Concini and C.\ Procesi in the early 1980s, yielding compactifications for certain symmetric varieties, and in particular, for semisimple algebraic groups of adjoint type \cite{DP}.  The use of these compactifications in the parametrization of the moduli space of bundles is one of the innovations of this paper.

The notion of a Hecke modification of a principal bundle is widely referred to in the literature (for example, see \cite{Norbury, KW, BZF}), however statements and results are often quite fragmented and given without justification, so a conscious attempt to systematize the exposition has been made in Section \ref{HMPB}.  After setting conventions with respect to root systems and weights, we discuss the loop group of a complex algebraic group and its corresponding affine Grassmannian, an infinite-dimensional homogeneous space.  The most relevant objects for us will be certain finite-dimensional subvarieties in the Grassmannian, known as Bruhat cells, which may be identified with certain double cosets in the loop group.  These Bruhat cells give the correct parameter spaces for Hecke modifications of a given bundle at a fixed point.  The structure theory here depends heavily on the work of Iwahori and Matsumoto \cite{IwahoriMatsumoto}, and to allow for a clear understanding of it, we review the notions of the affine root system and the affine Weyl group.  We then proceed to describe how these constructions can be made intrinsic to a point on a Riemann surface, and hence describe the spaces of Hecke modifications of a fixed principal bundle.  The section is concluded with the constructions of universal families of Hecke modifications of a fixed bundle, first for one and then for several modifications.

Section \ref{OWC} gives an overview of the construction of the wonderful compactification, largely following the treatment of S.\ Evens and B.F.\ Jones \cite{EJ}.  The structure of the ``standard'' open affine sets as well as the divisor at infinity are explicitly described.  The compactification admits left and right actions of the group analogous to those of $GL_n \C$ on the space of $n \times n$ matrices; we obtain explicit expressions for the associated infinitesimal actions, which become useful later for the deformation theory.  We also prove the existence of an involution extending the inversion map on $G$.

The main purpose that the wonderful compactification serves is in the development of the deformation theory for moduli of principal bundles, and this is carried out in Chapter \ref{UsingTheCompactifications}.  In the vector bundle case, when one bundle is given as a Hecke modification of another, there is still a map between the respective sheaves of sections. However, these sheaves of sections are not available to us in the principal bundle context, but what we can do is compactify the fibres of the fixed bundle, and consider families of bundles that map into this compactified bundle.  This construction allows us to define the sheaf whose global sections give us the infinitesimal deformations of the parameter space and to compute the Kodaira--Spencer map for the family of bundles we have constructed.

In Section \ref{Parametrization}, we give a general outline laying out sufficient conditions for when we obtain a parametrization of the moduli space.  The idea is that we introduce a number of modifications to the fixed bundle to obtain another bundle which is reducible to a maximal torus.  While this bundle will not be stable, if the family of bundles so constructed is of the right dimension, then nearby there will be an open set of stable bundles.  The surjectivity of the Kodaira--Spencer map amounts to the vanishing of the first cohomology of a certain vector bundle.  This vanishing requires that the locations of the modifications are chosen suitably generically.  This is also discussed in Section \ref{Parametrization}.

In the final section, we attempt to construct families satisfying the conditions developed previously in specific instances.  Unfortunately, because each Hecke modification introduces a certain number of parameters dependent on the root system, we are not always able to construct families of the requisite dimension, but are only able to obtain parametrizations for bundles with structure groups corresponding to the root systems of type $A_3, C_l$, and $D_l$ (i.e., the groups $PGL_4 \C, PSp_{2l} \C, PSO_{2l}\C$), and these only when the genus is even.  %This does not mean that the construction is of no use in general; indeed, we may always construct families which give submersions onto the moduli spaces, and for stack-theoretic purposes this is often sufficient.

One of the motivations for this paper was to extend results of I.\ Krichever \cite{Krichever} and Hurtubise \cite{Hurtubise}, which give a hamiltonian interpretation to the difference of two isomonodromic splittings on the moduli space of local systems, to the principal bundle case.  So as to maintain a reasonable length here, these considerations will be the subject of a forthcoming paper.

This paper is adapted from part of a doctoral thesis written under the supervision of Jacques Hurtubise.  I would like to thank him for his many ideas and for his encouragement over the course of innumerable discussions; these have contributed greatly to what appears here.

\section{Hecke Modifications of Principal Bundles} \label{HMPB}

\subsection{Notation for Roots and Weights} \label{rootsandweights}

Let $G$ be a semisimple algebraic group of rank $l$ over $\C$ and let $T \subseteq G$ be a maximal torus, $\g, \t$ their respective Lie algebras.  Let $\Phi$ be the corresponding root system and $W$ the associated Weyl group.  We will think of a root $\alpha \in \Phi$ as being either an element of the character group $X(T)$ or an element of $\t^*$ as the context dictates.  We will denote the root space corresponding to $\alpha$ by $\g_\alpha$.  The root lattice $\Lambda_r$ in the group of characters $X(T)$ will be denoted by $\Lambda_r$ and the weight lattice by $\Lambda$ and weights by $\lambda \in \Lambda$.  The following relation among these lattices holds:
\begin{align} \label{lattices}
\Lambda_r \subseteq X(T) \subseteq \Lambda \subseteq \t^*.
\end{align}
Coroots and coweights will be denoted by $\alpha^\vee$ and $\lambda^\vee$, respectively.  If $Y(T)$ is the group of cocharacters, then we have the dualization of (\ref{lattices}):
\begin{align} \label{colattices}
\Lambda_r^\vee \subseteq Y(T) \subseteq \Lambda^\vee \subseteq \t.
\end{align}
To be clear, if $\lambda \in X(T)$ and $\lambda^\vee \in Y(T)$ are thought of as homomorphisms $T \to \C^\times, \C^\times \to T$, respectively, then $\lambda \circ \lambda^\vee : \C^\times \to \C^\times$ is the map
\begin{align*}
z \mapsto z^{\langle \lambda^\vee, \lambda \rangle},
\end{align*}
where the pairing on the right side is the one we use when thinking of $\lambda$ and $\lambda^\vee$ as elements of $\t^*$ and $\t$, respectively.

The quotient $Y(T)/ \Lambda_r^\vee$ is called the fundamental group of $G$ and indeed it coincides with the topological fundamental group $\pi_1(G)$ \cite[Proposition 3.11.1]{DK}, which will be a finite abelian group.  The identification is obtained by restricting a cocharacter to $S^1 \subseteq \C^\times$ and taking the homotopy class.  Observe that this implies that the coroot lattice is precisely the subgroup of null-homotopic cocharacters.

A choice of a Borel subgroup $B$ containing $T$ (say with Lie algebra $\b$) is equivalent to a choice of a set of simple roots $\Delta := \{ \alpha_1, \ldots, \alpha_l \}$.  Let $\Phi^+, \Phi^-$ denote the corresponding sets of positive and negative roots, respectively.  Then there is a basis $\Delta^\vee := \{ \alpha_1^\vee, \ldots, \alpha_l^\vee \}$ of $\t$ such that in the natural pairing $\langle \, , \rangle : \t \otimes \t^* \to \C$, if
\begin{align*}
a_{ij} := \langle \alpha_i^\vee, \alpha_j \rangle
\end{align*}
then $A = (a_{ij})$ is the Cartan matrix (of finite type) from which $\g$ arises.\footnote{This is the convention taken in \cite{Kac}; one should note that the convention in \cite{Humphreys} is to take the transpose of this matrix.}  The set $\Delta^\vee$ gives a set of simple roots for the dual root system $\Phi^\vee \subseteq \t$.  The fundamental weights $\{ \lambda_i \}_{i=1}^l$ and coweights $\{ \lambda_i^\vee \}_{i=1}^l$ are bases of $\t^*$ and $\t$, respectively, dual to $\Delta$ and $\Delta^\vee$.  A coweight $\lambda^\vee \in \Lambda^\vee$ is called dominant if $\langle \lambda^\vee, \alpha_i \rangle \geq 0$ for all simple roots $\alpha_i, 1 \leq i \leq l$; clearly, this holds if and only if $\langle \lambda^\vee, \alpha \rangle \geq 0$ for all $\alpha \in \Phi^+$.  We will write $\Lambda_+, Y(T)_+$ and $\Lambda_{r+}$ for the sets of dominant weights, cocharacters and elements of the coroot lattice, respectively.  As is standard, we will denote by $\rho$ the half sum of the positive roots:  $2\rho = \sum_{\alpha \in \Phi^+} \alpha$; this coincides with the sum of the fundamental weights.

\subsection{Loop Groups and the Affine Grassmannian} 

\subsubsection{Definitions}

\paragraph{\emph{The Loop Group}}

Standard definitions and results on loop groups can be found in the book of Pressley and Segal \cite{PS}, which gives an analytic exposition, or in the work of Faltings \cite{Faltings2003} for an algebro-geometric one.  However, we will consider the following version as it is most amenable to our intended applications.  We will say that a map from an open subset of $\C$ to $G$ is \emph{meromorphic} if upon choosing an embedding $G \hookrightarrow GL_n \C$ the component functions are meromorphic on the open set.  Since these component functions for any such representation generate the coordinate ring of $G$, it is not hard to see that this is well-defined (in the sense that if the component functions are meromorphic in one representation, then they cannot acquire essential singularities in another).  On the other hand, the order of a pole for such a function is not a well-defined notion.  We will define the \emph{loop group} $L_\mero G = LG$ to be the group of germs of meromorphic $G$-valued functions at $0 \in \C$, the operation being pointwise multiplication in $G$; its elements will be called \emph{loops}.  The subgroup $L^+G \subseteq LG$ will be defined to be the subgroup of germs of holomorphic $G$-valued functions at $0$ and its elements will be referred to as \emph{positive loops}.  Observe that we may realize $G$ as a subgroup of $L^+G$ by considering the constant loops.  The cocharacter group $Y(T)$ may also be realized as a subgroup of $LG$.

If $K$ denotes the field of germs of meromorphic functions at $0$, and $R$ the ring of germs of holomorphic functions at $0$, then it is clear that $LG = G(K)$ is the set of $K$-valued points of $G$ and $L^+ G = G(R)$ is the set of $R$-valued points.  Fixing the standard coordinate $z$ on $\C$, we will typically identify $K$ with the field $\Czz$ of convergent Laurent series and $R$ with the convergent power series ring $\Cz$.  %As $R$ is a discrete valuation ring with quotient field $K$, we may apply the results of Iwahori and Matsumoto \cite[Section 2]{IwahoriMatsumoto} to $LG$ and $L^+G$.

Since a loop $\gamma \in LG$ is defined as a germ of a meromorphic function, its domain can always be taken to be a punctured disc centred at the origin.  As such, it defines a class in $\pi_1(G)$.  Clearly, elements of $L^+G$ define null-homotopic paths.

\begin{prop} \label{LGcomponents} \cite[Proposition 1.13.2 and its proof]{DK}
The map $LG \to \pi_1(G)$ which sends a loop to its homotopy class defines a group homomorphism whose kernel contains $L^+G$.  The connected components of $LG$ are indexed by $\pi_1(G)$.
\end{prop}

\

\paragraph{\emph{The Affine Grassmannian}} The \emph{loop} or \emph{affine Grassmannian} is defined as the homogeneous space
\begin{align*}
\Gr_G := LG/ L^+G,
\end{align*}
where we are simply quotienting by right multiplication.  $\Gr_G$ has the structure of a projective ind-variety, which means that there are projective varieties $X_j, j \in \N$ and closed immersions $X_j \hookrightarrow X_{j+1}$ such that $\Gr_G = \bigcup_{j \in \N} X_j$.  We will not go into how the ind-variety structure is defined (the interested reader may consult \cite[\S7.1]{Kumar}), but will later give descriptions of certain open sets in some of the subvarieties of $\Gr_G$.

By Proposition \ref{LGcomponents}, any two representatives of a class in $\Gr_G$ define the same homotopy class, which fact allows the following.

\begin{cor}
The connected components of $\Gr_G$ are indexed by $\pi_1(G)$.
\end{cor}

If $\varepsilon \in \pi_1(G)$, then $\Gr_G(\varepsilon)$ will denote the component of $\Gr_G$ corresponding to $\varepsilon$.  If $\varepsilon_1, \varepsilon_2 \in \pi_1(G)$, then there is a bijection $\Gr_G(\varepsilon_1) \arsim \Gr_G(\varepsilon_2)$ given by
\begin{align} \label{componentmap}
[ \sigma ] \mapsto [ \sigma \tilde{\varepsilon}_1^{-1} \tilde{\varepsilon}_2 ],
\end{align}
where $\tilde{\varepsilon}_i \in LG$ is a loop representing the homotopy class $\varepsilon_i$ for $i = 1,2$.

\subsubsection{The Affine Weyl Group and Affine Roots}

To a root system $\Phi$ with Weyl group $W$, we can associate an affine root system $\Phi_\af$ and affine Weyl group $W_\af$, which play the roles in the Bruhat decomposition of $LG$ that $W$ and $\Phi$ do in the decomposition of $G$.

\

\paragraph{\emph{The Affine Weyl Group}}

Recall that the Weyl group $W$ acts on the coweight lattice $\Lambda^\vee$ mapping the coroot lattice $\Lambda_r^\vee$ to itself.  We define the affine Weyl group to be the semi-direct product $W_\af := \Lambda_r^\vee \rtimes W$.  If $\lambda^\vee \in \Lambda_r^\vee$, we will often write $t(\lambda^\vee)$ when we think of it as an element of $W_\af$.  It is straightforward to check that
\begin{align} \label{corootconjugation}
w t(\lambda^\vee) w^{-1} = t(w \cdot \lambda^\vee).
\end{align}

\

\paragraph{\emph{The Affine Root System}}
The set of affine roots associated to $\Phi$ can be defined as $\Phi_\af := \Phi \times \Z$.  There is a decomposition $\Phi_\af = \Phi_\af^+ \coprod \Phi_\af^-$ into positive and negative roots, where
\begin{align*}
\Phi_\af^+ & := \Phi \times \Z_{> 0} \cup \Phi^+ \times \{ 0 \}, & \Phi_\af^- & := \Phi \times \Z_{< 0} \cup \Phi^- \times \{ 0 \}.
\end{align*}

The set $\Phi_\af$ carries an action of the group $W_\af$ which can be described as follows:  if $w \in W, \lambda^\vee \in \Lambda_r^\vee$, then
\begin{align} \label{affineWeylgroupaction}
w \cdot (\alpha, n) & = ( w \cdot \alpha, n), & t(\lambda^\vee) \cdot (\alpha, n) = (\alpha, n + \langle \lambda^\vee, \alpha \rangle).
\end{align}

\

\paragraph{\emph{Simple Reflections}}
If $\Phi = \Phi_1 \cup \cdots \cup \Phi_m$ is the decomposition of $\Phi$ into irreducible root systems $\Phi_j, 1 \leq j \leq m$, let $\theta_j$ denote the highest root in $\Phi_j$, and $\theta_j^\vee$ the corresponding coroot.  With this notation, $W_\af$ is a Coxeter group with involutive generators
\begin{align*}
S := \{ s_{0,1}, \ldots, s_{0,m}, s_1, \ldots, s_l \},
\end{align*}
where $s_i, 1 \leq i \leq l$ are the simple reflections corresponding to the simple roots (the usual generators for $W$), $s_{\theta_j}, 1 \leq j \leq m$ are the reflections corresponding to the roots $\theta_j$, and
\begin{align*}
s_{0,j} := s_{\theta_j} t(\theta_j^\vee).
\end{align*}
Since the coroot lattice is generated by the $\Z$-span of the $W$-orbits of the $\theta_j^\vee, 1 \leq j \leq m$, by (\ref{corootconjugation}) these do indeed generate $W_\af$.  We will set $\alpha_{0,j} := ( - \theta_j, 1) \in \Phi_\af^+$.

\

\paragraph{\emph{The Length Function}}
There is a length function $\ell : W_\af \to \N$ which takes $s \in W_\af$ to the smallest $k$ such that $s$ can be written as a product of $k$ elements of $S$.  This extends the usual length function on $W$.  If $\ell(s) = k$ and $s = s_{i_1} \ldots s_{i_k}$ is an expression with $i_j \in \{ (0,1), \ldots, (0,m), 1, \ldots, l \}$, then this expression is called reduced.  For $\sigma \in W_\af$, let us denote
\begin{align*}
\Phi_\af^s := \{ \beta \in \Phi_\af^+ \, : \, s^{-1} \beta \in \Phi_\af^- \}.
\end{align*}

\begin{lem} \label{AWGlem}
    \begin{enumerate}
    \item[(a)] For $i \in I$, $\Phi_\af^{s_i} = \{ \alpha_i \}$.
    \item[(b)] \cite[Lemma 1.3.14]{Kumar} If $s \in W_\af$, and $s = s_{i_1} \ldots s_{i_k}$ is a reduced expression, then
    \begin{align*}
    \Phi_\af^s = \{ \alpha_{i_1}, s_{i_1} \alpha_{i_2}, s_{i_1} s_{i_2} \alpha_{i_3}, \ldots, s_{i_1} \ldots s_{i_{k-1}} \alpha_{i_k} \}.
    \end{align*}
    In particular,
    \begin{align*}
    \ell(s) = \# \Phi_\af^s.
    \end{align*}
    \item[(c)] If $s \in W_\af$ and $i \in I$, then $\alpha_i$ lies in exactly one of $\Phi_\af^s$ or $\Phi_\af^{s_i s}$, and $\alpha_i \in \Phi_\af^s$ if and only if $\ell(s_i s) < \ell(s)$.
    \item[(d)] Let $W_s := \{ w \in W \, : \, ws = sv$ for some $v \in W \}$.  Then for $\lambda^\vee \in \Lambda_r^\vee$,
    \begin{align*}
    W_{t(\lambda^\vee)} = \{ w \in W \, : \, w t(\lambda^\vee) = t(\lambda^\vee) w \} = \{ w \in W \, : \, w \cdot \lambda^\vee = \lambda^\vee \}.
    \end{align*}
    \item[(e)] If $\lambda^\vee \in \Lambda_{r+}$ and $s_i \not\in W_{t(\lambda^\vee)}$, then $\ell(s_i t(\lambda^\vee)) < \ell(t(\lambda^\vee))$.  More generally,
        \begin{align*}
        \ell\big( t(\lambda^\vee) \big) = \max \{ \ell\big( w t(\lambda^\vee) \big) \, : \, w \in \lfloor W/ W_{t(\lambda^\vee)} \rfloor \},
        \end{align*}
        if $\lfloor W/ W_{t(\lambda^\vee)} \rfloor$ is a set of coset representatives of minimal length.
    \end{enumerate}
\end{lem}

\begin{proof}
If $1 \leq i \leq l$, then if $(\alpha, n) \in \Phi_\af^+$ and $s_i (\alpha, n) = (s_i \alpha, n) \in \Phi_\af^-$, it follows that $n = 0$ and $s_i \alpha \in \Phi^-$, so $\alpha = \alpha_i$.  By dealing with each irreducible component separately, we may assume that $\Phi$ is irreducible and that the remaining simple root is $\alpha_0$.  Suppose $(\alpha, n) \in \Phi_\af^+$ and $s_0 (\alpha, n) = s_\theta t(\theta^\vee) (\alpha, n) = (s_\theta \alpha, n + \langle \alpha, \theta^\vee \rangle) \in \Phi_\af^-$.  Since $\theta$ is a long root $\langle \alpha, \theta^\vee \rangle \in \{ 0, \pm 1 \}$; also, it is clear that $\langle \alpha, \theta^\vee \rangle \leq 0$, so in fact, $\langle \alpha, \theta^\vee \rangle \in \{ 0, 1 \}$.  If $\langle \alpha, \theta^\vee \rangle = 0$, then $s_\theta \alpha = \alpha$ and $n = 0$, so we get $s_0 (\alpha, 0) = (\alpha, 0) \in \Phi_\af^+$, a contradiction.  If $\langle \alpha, \theta^\vee \rangle = -1$, then $s_\theta \alpha = \alpha + \theta \in \Phi$ forces $\alpha \in \Phi^-$ and so $n > 0$, but then $(\alpha + \theta, n + \langle \alpha, \theta^\vee \rangle ) \in \Phi_\af^+$, again a contradiction.  It follows that $\alpha = \pm \theta$, and since $\langle \alpha, \theta^\vee \rangle \leq 0$, we must have $\alpha = -\theta$, in which case $\langle \alpha, \theta^\vee \rangle =-2$.  This means that $(s_\theta \alpha, n + \langle \alpha, \theta^\vee \rangle ) = (\theta, n-2)$, and the only possibility is $n = 1$.  Hence $(\alpha, n) = (-\theta, 1) = \alpha_0$.  Part (b) is obtained from (a) by a straightforward induction.

Observe that $(s_i s)^{-1} \alpha_i = s^{-1} s_i \alpha_i = - s^{-1} \alpha_i$.  This implies that $\alpha_i$ lies in precisely one of $\Phi_\af^s$ or $\Phi_\af^{s_i s}$.  If $\alpha_i \in \Phi_\af^s$, then one can check that
\begin{align*}
\beta \mapsto s_i \beta
\end{align*}
gives an injection $\Phi_\af^{s_i s} \to \Phi_\af^s$, but $\Phi_\af^s$ contains one more element.  The converse is exactly the same.  This proves (c).

Clearly, if $w t(\lambda^\vee) = t(\lambda^\vee)w$, then $w \in W_{t(\lambda^\vee)}$.  Conversely, if $w t(\lambda^\vee) = t(\lambda^\vee) v$ for some $v \in W$, then writing $w t(\lambda^\vee) = t( w \cdot \lambda^\vee) w = t(\lambda^\vee) v$, by uniqueness of the factorization in a semi-direct product, it follows that $w \cdot \lambda^\vee = \lambda^\vee$ and $v = w$.

If $s_i \not\in W_{t(\lambda^\vee)}$, then since $\lambda^\vee \neq s_i \lambda^\vee = \lambda^\vee - \langle \lambda^\vee, \alpha_i \rangle \alpha_i^\vee$, it follows that $\langle \lambda^\vee, \alpha_i \rangle > 0$, since $\lambda^\vee \in \Lambda_{r+}$.  But then, $t(\lambda^\vee)^{-1} \cdot (\alpha_i, 0) = (\alpha_i, - \langle \lambda^\vee, \alpha_i \rangle) \in \Phi_\af^-$ and so $\alpha_i \in \Phi_\af^{t(\lambda^\vee)}$.  So the first part of (e) follows from (c).

Suppose $w \in \lfloor W/ W_{t(\lambda^\vee)} \rfloor$.  Then we can write $w = s_i v$ for some $1 \leq i \leq l$ and $v \in \lfloor W/ W_{t(\lambda^\vee)} \rfloor$ with $\ell(v) = \ell(w)-1$.  By (c), $\alpha_i \in \Phi_\af^w$, i.e.\ $w^{-1} \alpha_i \in \Phi^-$.  Then $t(\lambda^\vee)^{-1} w^{-1} \alpha_i \in \Phi_\af^-$, hence $\alpha_i \not\in \Phi_\af^{w t(\lambda^\vee)}$.  Therefore $\alpha_i \in \Phi_\af^{v  t(\lambda^\vee)}$, and so again by (c), $\ell( wt(\lambda^\vee)) < \ell( v t(\lambda^\vee) ) \leq \ell( t(\lambda^\vee) )$, by induction, with equality if and only if $v = e$.
\end{proof}

\

\paragraph{\emph{Interpretation In Terms of Affine Transformations}} The affine Weyl group is commonly described as a group of affine transformations of the real vector space $\t_\R \subseteq \t$ spanned by $\Phi^\vee$.  In this realization, $W$ acts in the usual manner and $\Lambda_r^\vee$ acts by translations.  If $\alpha \in \Phi, n \in \Z$, let $P_{\alpha, n}$ denote the hyperplane
\begin{align*}
P_{\alpha, n} := \{ \lambda^\vee \in \t_\R \, : \, \langle \lambda^\vee, \alpha \rangle = n \}.
\end{align*}
With this notation, the elements $s_{0,j}, 1 \leq j \leq m$ correspond to reflections in the planes $P_{\alpha_{0,j}, 1}$.  Since $P_{-\alpha, -n} = P_{\alpha, n}$, we may always assume that $\alpha \in \Phi^+$.  Any plane $P_{\alpha, n}$ divides $\t_\R$ into two half-planes
\begin{align*}
P_{\alpha, n}^+ & := \{ \lambda^\vee \in \t_\R \, : \, \langle \lambda^\vee, \alpha \rangle \geq n \}, & P_{\alpha, n}^- & := \{ \lambda^\vee \in \t_\R \, : \, \langle \lambda^\vee, \alpha \rangle \leq n \}.
\end{align*}
With this, it is clear that $W_\af$ acts on the set $\Phi_\af'$ of such half-planes.  It is not hard to see that $\Phi_\af'$ and the set $\Phi_\af$ of affine roots are isomorphic as $W_\af$-sets.

A Weyl alcove is defined as a connected component of the complement of $\bigcup_{\alpha \in \Phi, n \in \Z} P_{\alpha, n}$.  It is a basic fact, though one that we will not need, that $W_\af$ acts simply transitively on the set of the Weyl alcoves \cite[Corollary 1.8]{IwahoriMatsumoto}.
%The fundamental Weyl alcove is defined as the intersection
%\begin{align*}
%P_{\alpha_{0,1},1}^- \cap \cdots \cap P_{\alpha_{0,m},1} \cap P_{\alpha_1,0}^+ \cap P_{\alpha_2,0}^+ \cap \cdots \cap P_{\alpha_l,0}^+.
%\end{align*}
%

\subsubsection{Bruhat Decomposition}

As a set of left cosets of $L^+G$ in $LG$, it is clear that $\Gr_G$ admits a left $L^+ G$-action.  Understanding the $L^+G$ orbits in $\Gr_G$ thus amounts to understanding the double $L^+G$ cosets in $LG$. The orbits we are particularly interested in are those of the dominant cocharacters $\lambda^\vee \in Y(T)$; such orbits will be denoted $\Gr_G^{\lambda^\vee}$ and are called the \emph{Bruhat cells} of $\Gr_G$. Since $LG$ is a group over the field $K$ with a discrete valuation, we may apply the results of \cite{IwahoriMatsumoto} to obtain a Bruhat decomposition.  To this end, we now introduce some relevant subgroups of $LG$.

\

\paragraph{\emph{The Extended Affine Weyl Group}}

If $N$ is the normalizer of $LT$ in $LG$, then the extended affine Weyl group is defined as
\begin{align*}
\widetilde{W}_\af := N/L^+T \cong Y(T) \rtimes W.
\end{align*}
Going back to (\ref{colattices}), we see that $\widetilde{W}_\af$ contains the affine Weyl group $W_\af = \Lambda_r^\vee \rtimes W$ as a subgroup.  Since for any $\lambda^\vee \in \Lambda^\vee, w \in W$, $\lambda^\vee - w \cdot \lambda^\vee \in \Lambda_r^\vee$, it follows that $W_\af$ is normal in $\widetilde{W}_\af$ with quotient isomorphic to $Y(T)/\Lambda_r^\vee \cong \pi_1(G)$.  In fact, there is a subgroup of $Y(T)$ which maps isomorphically onto $\pi_1(G)$ under the above quotient so that $\widetilde{W}_\af \cong W_\af \rtimes \pi_1(G)$.  This subgroup may be realized as the elements of $\widetilde{W}_\af$ which map the fundamental Weyl alcove to itself \cite[\S1.7]{IwahoriMatsumoto}.  Since any element of $\widetilde{W}_\af$ can be written uniquely in the form $s \varepsilon$ with $s \in W_\af, \varepsilon \in \pi_1(G)$, we can extend the length function $\ell$ to $\widetilde{W}_\af$ by setting
\begin{align*}
\ell(s \varepsilon) = \ell(s).
\end{align*}

\

\paragraph{\emph{Root Groups}}

Given a root $\alpha \in \Phi$, consider the root groups $U_\alpha \subseteq G$, where $U_\alpha = \exp ( \C \xi_\alpha)$ for a root vector $\xi_\alpha \in \g_\alpha$.  Given $n \in \Z$, we may restrict the isomorphism
\begin{align*}
\mathbb{G}_a (K) \xrightarrow{\exp} LU_\alpha = U_\alpha(K) \subseteq LG,
\end{align*}
to the additive subgroup $\C z^n \subseteq \mathbb{G}_a(K)$.  The image will be a subgroup of $LU_\alpha$ isomorphic to $\mathbb{G}_a(\C) = \C$.  We will denote this subgroup by
\begin{align*}
U_{\alpha, n} = \exp( \C z^n \xi_\alpha),
\end{align*}
so that the set of all such subgroups is indexed by elements of $\Phi_\af$.  If $w \in W$ and $\widetilde{w} \in N$ is a representative, then
\begin{align*}
\Ad \, \widetilde{w} (U_{\alpha, n}) = U_{w \cdot \alpha, n},
\end{align*}
and if $\lambda^\vee \in Y(T)$ then
\begin{align*}
\Ad \, \lambda^\vee (U_{\alpha, n}) = U_{\alpha, n + \langle \lambda^\vee, \alpha \rangle}.
\end{align*}
Thus, $\widetilde{W}_\af$ permutes these root groups in such a way that the restriction to $W_\af$ acts on the indices as in (\ref{affineWeylgroupaction}).

There is an evaluation map $\ev : L^+ G \to G$ which simply takes a germ to its value at $0$. We will denote by $I := \ev^{-1}(B)$ its pre-image in $L^+G$; this is the Iwahori subgroup.  Using the arguments of \cite[Corollary 2.7(ii),(iii)]{IwahoriMatsumoto}, one can choose coset representatives for $I$ as follows.

\begin{prop} \label{simpleIcoset}
    \begin{enumerate}
    \item[(a)] For $1 \leq i \leq l$, we have $I s_i I = U_{\alpha_i,0} s_i I$.
    \item[(b)] For $1 \leq j \leq m$, we have $I s_{0,j} I = U_{-\theta_j, 1} s_{0,j} I$.
    \end{enumerate}
Therefore, if we identify $\alpha \in \Phi$ with $(\alpha,0) \in \Phi_\af$ and $(-\theta_j,1)$ with $\alpha_{0,j}$, then if $i$ lies in the index set $\{ (0,1), \ldots, (0,m), 1, \ldots, l \}$, we have
\begin{align*}
I s_i I = U_{\alpha_i} s_i I.
\end{align*}
\end{prop}

\

\paragraph{\emph{Bruhat Decomposition}}

We will use the following information about double cosets in $LG$ in our main result about the $\Gr_G^{\lambda^\vee}$.

\begin{prop} \label{cosetdecomposition}
Let $s \in W_\af$.  Then
   \begin{align*}
   Is I = \left( \prod_{\beta \in \Phi_\af^s } U_\beta \right) s I,
   \end{align*}
   and
   \begin{align*}
   L^+G s L^+G = \coprod_{w \in \lfloor W/W_s \rfloor } \left( \prod_{\beta \in \Phi_\af^{w s}} U_\beta \right) w s L^+ G,
   \end{align*}
   where $W_s$ is as in Lemma \ref{AWGlem}(d).
\end{prop}

\begin{proof}
Suppose $s = s_{i_1} \cdots s_{i_k}$ is a reduced expression.  Then using an argument of \cite[\S\S3,5]{Lansky}, by repeated application of Proposition \ref{simpleIcoset}, we obtain
\begin{align*}
I s I & = I s_{i_1} \cdots s_{i_k} I = I s_{i_1} I s_{i_2} I \cdots I s_{i_k} I = U_{\alpha_{i_1}} s_{i_1} I s_{i_2} I \cdots I s_{i_k} I = U_{\alpha_{i_1}} s_{i_1} U_{\alpha_{i_2}} s_{i_2} I \cdots I s_{i_k} I \\
& = U_{\alpha_{i_1}} s_{i_1} U_{\alpha_{i_2}} s_{i_2} \cdots U_{\alpha_{i_k}} s_{i_k} I = U_{\alpha_{i_1}} U_{s_{i_1} \alpha_{i_2}} \cdots U_{s_{i_1} \cdots s_{i_{k-1}} \alpha_{i_k}} s I = \left( \prod_{\beta \in \Phi_\af^s } U_\beta \right) s I.
\end{align*}
This is the first equality.

For the second, we will begin by noting that it follows from \cite[Proposition 2.4]{IwahoriMatsumoto} that $L^+ G = IWI$.  Hence,
\begin{align*}
L^+G s L^+G = \coprod_{w \in W} I w I s L^+G = \bigcup_{w \in W} I w s I L^+ G = \bigcup_{w \in W} I w s L^+ G.
\end{align*}
Observe that $W_s$ is precisely the set of $w \in W$ for which $I w s L^+ G = I s L^+G$.  So the above union need only be taken over a set of coset representatives, which we may assume to be minimal in their respective cosets.  Hence,
\begin{align*}
L^+G s L^+G = \coprod_{w \in \lfloor W/W_s \rfloor } I w s L^+ G = \coprod_{w \in \lfloor W/W_s \rfloor } \left( \prod_{\beta \in \Phi_\af^{w s}} U_\beta \right) w s L^+ G.
\end{align*}
\end{proof}

We now come to some of the most relevant information about the $\Gr_G^{\lambda^\vee}$ for us.

\begin{thm} \label{Bruhatdecomp}
    \begin{enumerate}
    \item[(a)] The affine Grassmannian is a disjoint union of the $\Gr_G^{\lambda^\vee}$:
    \begin{align*}
    \Gr_G = \coprod_{\lambda^\vee \in Y(T)_+ } \Gr_G^{\lambda^\vee}.
    \end{align*}
    \item[(b)] $\Gr_G^{\lambda^\vee}$ is a rational variety of dimension
    \begin{align*}
    \dim \Gr_G^{\lambda^\vee} = \ell\big( t(\lambda^\vee) \big) = \sum_{\alpha \in \Phi^+} \langle \lambda^\vee, \alpha \rangle = 2 \langle \lambda^\vee, \rho \rangle.
    \end{align*}
    In fact,
    \begin{align*}
    V := \left( \prod_{\beta \in \Phi_\af^{t(\lambda^\vee)}} U_\beta \right) \lambda^\vee L^+G
    \end{align*}
    is an open set in $\Gr_G^{\lambda^\vee}$ and
    \begin{align} \label{openaffinecover}
    \Gr_G^{\lambda^\vee} = \bigcup_{w \in \lfloor W/W_{t(\lambda^\vee)} \rfloor } V_w,
    \end{align}
    where $V_w := w \cdot V$ is the $w$-translation of $V$, is an open covering by affine spaces.
    \item[(c)] Each $\Gr_G^{\lambda^\vee}$ is a homogeneous space for a group $G( \C[z]/(z^{n+1}))$ for some $n \geq 0$.
    \end{enumerate}
\end{thm}

\begin{proof}
The statement in (a) follows directly from \cite[Corollary 2.35(ii)]{IwahoriMatsumoto} (cf.\ \cite[Proposition 8.6.5]{PS}).

For (b), first we note that given $\lambda^\vee \in Y(T)_+$, $\Gr_G^{\lambda^\vee} = \Gr_G^{\lambda_1^\vee \varepsilon}$ for some $\lambda_1^\vee \in \Lambda_r^\vee, \varepsilon \in \pi_1(G)$.  But then we will have an isomorphism $\Gr_G^{\lambda_1^\vee} \cong \Gr_G^{\lambda_1^\vee \varepsilon}$ via a map as in (\ref{componentmap}), so it suffices to assume that $\lambda^\vee \in \Lambda_{r+} \subseteq W_\af$.  In this case, it follows from Proposition \ref{cosetdecomposition} that $\Gr_G^{\lambda^\vee}$ is a finite union of affine spaces.  The largest one will be an open set and its dimension will give the dimension of $\Gr_G^{\lambda^\vee}$; but by Lemma \ref{AWGlem}(e), this is precisely the set $V$.
It is straightforward to check that
\begin{align*}
\Phi_\af^{t(\lambda^\vee)} = \{ (\alpha, n) \in \Phi^+ \times \Z_{\geq 0} \, : \, 0 \leq n \leq \langle \lambda^\vee, \alpha \rangle -1 \},
\end{align*}
and so
\begin{align*}
\dim \Gr_G^{\lambda^\vee} = \# \Phi_\af^{t(\lambda^\vee)} = \ell\big( t(\lambda^\vee) \big) = \sum_{\alpha \in \Phi^+} \langle \lambda^\vee, \alpha \rangle = 2 \langle \lambda^\vee, \rho \rangle.
\end{align*}

Observe that for $w \in \lfloor W/W_{t(\lambda^\vee)} \rfloor$, we have $w^{-1} \phi_\af^{w t(\lambda^\vee)} \subseteq \Phi_\af^{t(\lambda^\vee)}$, so that
\begin{align*}
\left( \prod_{\beta \in \Phi_\af^{w t(\lambda^\vee)}} U_\beta \right) w \lambda^\vee L^+ G = w \cdot \left( \prod_{\beta \in \Phi_\af^{w t(\lambda^\vee)}} U_{w^{-1} \beta} \right) \lambda^\vee L^+ G  \subseteq w \cdot V,
\end{align*}
and hence the expression in (\ref{openaffinecover}) then comes from Proposition \ref{cosetdecomposition}.

If we let $R_n = \C[z]/(z^{n+1})$, then there are natural maps $R \to R_n$ for $n \geq 0$, which yield natural homomorphisms $\pi_n : L^+ G = G(R) \to G( R_n)$.  Then (c) amounts to saying that there is some $n$ such that $\ker \pi_n$ lies in the isotropy of $\lambda^\vee$.  But for a fixed $\lambda^\vee$, there will be a finite number of root groups $U_\beta$ that appear in a double coset decomposition as in Proposition \ref{cosetdecomposition}.  One sees that taking $n$ large enough, an element of $\ker \pi_n$, after rearranging factors, will not have any component in these root groups, so it does indeed stabilize $\lambda^\vee$.
\end{proof}

\begin{rmk}
An alternative way of computing $\dim \Gr_G^{\lambda^\vee}$ is indicated in \cite[\S2.2]{Norbury}.  The isotropy group of $\lambda^\vee$ in $L^+ G$ is readily computed to be $L^+ G \cap (\Ad \, \lambda^\vee ) L^+ G$, and so the tangent space to $\Gr_G^{\lambda^\vee}$ at $[\lambda^\vee]$ can be identified with $L^+ \g / \big( L^+ \g \cap (\Ad \, \lambda^\vee ) L^+ \g \big)$.  But now the computation of the dimension of this space is essentially the same as finding $\ell(t(\lambda^\vee))$.
\end{rmk}

The proof of the theorem shows us a way to choose coset representatives for $\Gr_G^{\lambda^\vee}$ over an open set isomorphic to an affine space.  Since it is a homogeneous space, this open set can be translated to obtain an open covering of $\Gr_G^{\lambda^\vee}$.  Therefore we may record the following.

\begin{cor} \label{LGsection}
The projection maps $\pi : L^+G \cdot \lambda^\vee \cdot L^+G \to \Gr_G^{\lambda^\vee}$ admit local sections:  for $w \in \lfloor W/ W_{t(\lambda^\vee)} \rfloor$, there are $f_w : V_w \to L^+G \cdot \lambda^\vee \cdot L^+G$ such that $\pi \circ f_w = \1_{V_w}$ and for which the loop $f_w(\sigma)$ is convergent on $\C^\times$, for all $\sigma \in V_w$.
\end{cor}

\begin{proof}
The statement about the convergence comes from noting that the value of $f_w(\sigma)$ is a finite product of elements of the root groups, each of which converges on all of $\C$, multiplied by a cocharacter, which is convergent on $\C^\times$.
\end{proof}

\begin{cor} \label{coordchange}
The group of local changes of the coordinate $z$ acts holomorphically on $\Gr_G^{\lambda^\vee}$
\end{cor}

\begin{proof}
By Theorem \ref{Bruhatdecomp}(c), the action factors through $\Aut \, \C[z]/(z^{n+1}) \cong (\C[z]/(z^{n+1}))^\times$, which acts algebraically on $\Gr_G^{\lambda^\vee}$.
\end{proof}

\subsubsection{Intrinsic Grassmannians}

Let $X$ be a Riemann surface and let $x \in X$.  Let $\G'(x)$ be the sheaf on $X$ whose value at $U \subseteq X$ is the group of holomorphic $G$-valued maps $U \setminus \{ x \} \to G$.  Let $\G(x) \subseteq \G'(x)$ be the sheaf of meromorphic $G$-valued functions with poles only at $x$ (we determine whether a $G$-valued function on $X$ is meromorphic by choosing a coordinate centred at $x$; clearly, this is independent of the choice of coordinate).  Let $\G$ be the subsheaf of holomorphic $G$-valued functions.  Then we may consider the stalks $\G(x)_x, \G_x$ at $x$ and the quotient
\begin{align*}
\Gr_G(x) := \G(x)_x / \G_x.
\end{align*}
Observe that $\G_x = G(\OO_{X,x}), \G(x)_x = G( \mathcal{K}_{X,x})$, where $\OO_{X,x}$ is the (analytic) local ring at $x$ and $\mathcal{K}_{X,x}$ its quotient field.  A choice of coordinate $z$ centred at $x$ fixes an isomorphism $\OO_{X,x} \arsim \C \{\!\{ z \}\!\}$ and hence isomorphisms
\begin{align*}
\G(x)_x & \xrightarrow{\sim} LG, & \G_x & \xrightarrow{\sim} L^+ G,
\end{align*}
finally yielding one
\begin{align} \label{intGmiso}
\Gr_G(x) \xrightarrow{\sim} LG / L^+ G = \Gr_G.
\end{align}
We would like to stratify these intrinsic Grassmannians $\Gr_G(x)$ as in (\ref{Bruhatdecomp}) simply by transporting the stratification across one of these isomorphisms.  However, to legitimize this, we need to see that the types are independent of the choice of coordinate.

Suppose $\gamma : U \setminus \{ x \} \to G$ represents a class in $\Gr_G(x)$ (where $U$ is a neighbourhood of $x$); upon choosing a coordinate $z$, we may assume $\gamma(z) = \gamma_+(z) \lambda^\vee(z)$ for some $\lambda^\vee \in Y(T)$.  We observe that there is no canonical group structure on $U$ (or any of its subsets), so it does not make sense to think of $\lambda^\vee$ as a homomorphism unless a coordinate is chosen.  But once we do, and realize an isomorphism $T \cong (\C^\times)^l$, then we may write
\begin{align*}
\lambda^\vee(z) = ( z^{r_1}, \ldots, z^{r_l})
\end{align*}
for some $r = (r_1, \ldots, r_l) \in \Z^l$; indeed $\lambda^\vee$ is determined by $r$.  If $w$ is another choice of coordinate, then
\begin{align*}
z = z(w) = w f(w),
\end{align*}
for some holomorphic nowhere-vanishing function $f(w)$.  Then
\begin{align*}
\lambda^\vee \big( z(w) \big) = \lambda^\vee(w) \lambda^\vee \big( f(w) \big).
\end{align*}
We note then that $\gamma_+ ( w f(w) ), \lambda^\vee ( f(w) ) \in L^+ G$, and so
\begin{align*}
\gamma(w) = \gamma_+ \big( w f(w) \big) \lambda^\vee(w) \lambda^\vee \big( f(w) \big)
\end{align*}
is also of type $\lambda^\vee$.

\begin{prop} \label{intBruhatdecomp}
There are Bruhat decompositions
\begin{align*}
\Gr_G(x) = \coprod_{ \lambda^\vee \in \Lambda_+^\vee } \Gr_G^{\lambda^\vee}(x),
\end{align*}
where $\Gr_G^{\lambda^\vee}(x)$ is the preimage of $\Gr_G^{\lambda^\vee}$ under an isomorphism (\ref{intGmiso}).
\end{prop}

\subsection{Hecke Modifications} \label{HMPBmodfn}

We will let $G$ be as above and fix a principal $G$-bundle $Q$ over $X$.  P.\ Norbury \cite{Norbury} defines a \emph{Hecke modification} of a principal bundle $Q$ (supported) at $x \in X$ as a pair $(P, s)$ consisting of a $G$-bundle $P$ and an isomorphism
\begin{align} \label{Hmdefn}
s : P|_{X_0} \to Q|_{X_0},
\end{align}
where $X_0 := X \setminus \{ x \}$.  We will want to restrict this definition somewhat so as to make sense of a meromorphic modification.  Let $X_1$ be a neighbourhood of $x$ on which we can choose trivializations $\psi_1$ and $\varphi_1$ of $Q|_{X_1}$ and $P|_{X_1}$, respectively.  Since $s|_{X_0}$ is an isomorphism, the composition
\begin{align*}
X_{01} \times G \xrightarrow{ \varphi_1^{-1}} P|_{X_{01}} \xrightarrow{s} Q|_{X_{01}} \xrightarrow{\psi_1} X_{01} \times G
\end{align*}
is an isomorphism of trivial $G$-bundles over $X_{01}$, so is of the form $\1_{X_{01}} \times L_{\sigma}$ for some holomorphic $\sigma : X_{01} \to G$.  Shrinking $X_1$ simply restricts $\sigma$, and so $\sigma$ defines an element in the stalk $\G'(x)_x$, which we will also denote by $\sigma$.  We will say that the modification is \emph{meromorphic} if $\sigma \in \G(x)_x$, i.e., if $\sigma$ is a meromorphic germ.

It is not hard to see that a change in the trivialization for $P$ or $Q$, respectively, amounts to multiplying $\sigma$ on the right or the left, respectively, by a positive loop, so our definition of a modification as being meromorphic is independent of the choice of trivializations.  In all that follows, we will assume that we are working with meromorphic Hecke modifications.

By Proposition \ref{intBruhatdecomp}, $[\sigma] \in \Gr_G^{\lambda^\vee}(x)$ for some dominant coweight $\lambda^\vee$.  A change in the $P$-trivialization amounts to right multiplication, so does not change the class of $\sigma$ at all; changing the $Q$-trivialization is the same as multiplication on the left, so we remain in the $\G_x$-orbit $\Gr_G^{\lambda^\vee}(x)$.  Therefore the orbit coweight $\lambda^\vee$ is independent of the choices made, so we can define the modification to be of \emph{type} $\lambda^\vee$.  Indeed, once we have chosen a coordinate $z$, by choosing the trivializations $\varphi_1, \psi_1$ appropriately, we may assume $\sigma(z) = \lambda^\vee(z)$ is itself a cocharacter.  If $G$ is of adjoint type, then $Y(T) = \Lambda^\vee$, so the fundamental coweights are all cocharacters.  In this case, a modification will be called \emph{simple} if its type $\lambda^\vee$ is one of the fundamental coweights.

If we choose trivializations $\psi_0$ and $\varphi_0$ of $Q|_{X_0}$ and $P|_{X_0}$, respectively (this is possible by a theorem of Harder \cite{Harder}), then we can form the respective transition functions $h_{01}, g_{01}$.  These are related by
\begin{align} \label{tfreln}
g_{01} = h_{01} \sigma.
\end{align}

We will say that two modifications $s_1 : P_1 |_{X_0} \to Q|_{X_0}$ and $s_2 : P_2 |_{X_0} \to Q|_{X_0}$ are \emph{equivalent} or \emph{isomorphic} if there is an isomorphism $\alpha : P_1 \to P_2$ and a commutative diagram
\begin{align*}
\xymatrixrowsep{3pc}\xymatrixcolsep{1pc}
\xymatrix{
P_1|_{X_0} \ar[rr]^\alpha \ar[dr]_{s_1} & & P_2 |_{X_0} \ar[dl]^{s_2} \\
& Q|_{X_0}. & }
\end{align*}
The following statement is straightforward to prove.

\begin{lem}
Two Hecke modifications $(P_1, s_1), (P_2, s_2)$ of $Q$ are isomorphic if and only if the corresponding $\sigma_1, \sigma_2$ (using the same choices of trivializations) as constructed above yield the same class in $\Gr_G(x)$.
\end{lem}

\subsubsection{Topological Considerations}

Using the notation above, then the topological types of $P$ and $Q$, respectively are given by the homotopy classes of $g_{01}$ and $h_{01}$ \cite[proof of Proposition 5.1]{Ramanathanthesis}, but these are related by (\ref{tfreln}), so the following relationship arises as a result of Proposition \ref{LGcomponents}.

\begin{prop} \label{HMtoptype}
If $\varepsilon(P), \varepsilon(Q) \in \pi_1(G)$ represent the topological types of $P$ and $Q$, respectively, and $P$ is obtained from $Q$ by introducing a Hecke modification of type $\lambda^\vee$, then
\begin{align*}
\varepsilon(P) = \varepsilon(Q) + \varepsilon(\lambda^\vee),
\end{align*}
if the modification lies in $\Gr_G( \varepsilon(\lambda^\vee))$, i.e.\ $\lambda^\vee$ yields the class $\varepsilon(\lambda^\vee) \in \pi_1(G)$.
\end{prop}

Here we have used additive notation for $\pi_1(G)$.

\subsection{Spaces of Hecke Modifications}

\subsubsection{Meromorphic Sections and Modifications}

We define a section of $Q$ over an open $U \subseteq X$ to be \emph{meromorphic} if when composed with a trivialization $Q|_U \xrightarrow{\sim} U \times G$, and thus written in the form
\begin{align*}
\1_U \times \gamma
\end{align*}
for some $G$-valued function $\gamma : U \to G$, $\gamma$ is meromorphic.

Fix $x \in X$.  We may consider the sheaf of sets $\L_x Q$ whose value at an open $U \subseteq X$ is the set of meromorphic sections of $Q$ over $U$ with poles only at $x$.  Then there is a point-wise right action of $\G$ on $\L_x Q$.  Then we can identify the space of (meromorphic) Hecke modifications of $Q$ supported at $x$ with the quotient of the stalks
\begin{align*}
\Gr_Q(x) := (\L_x Q)_x / \G_x.
\end{align*}

\begin{prop} \label{Hmspace}
The set $\Gr_Q(x)$ corresponds precisely to the set of all (meromorphic) Hecke modifications of $Q$ supported at $x$.
\end{prop}

\begin{proof}
Given $\varsigma \in \Gr_Q(x)$, we can choose a representative section $\tilde{\varsigma}$ of $Q$ over some small disc $X_1$ centred at $x$.  Let $\psi_1$ be a trivialization of $Q$ over $X_1$ and consider the meromorphic map $\sigma : X_1 \to G$ defined by
\begin{align*}
X_1 \xrightarrow{ \tilde{\varsigma} } Q|_{X_1} \xrightarrow{ \psi_1 } X_1 \times G \xrightarrow{p_G} G.
\end{align*}
Then $\sigma$ is holomorphic on $X_{01}$ so choosing a trivialization $\psi_0$ of $Q$ on $X_0 = X \setminus \{ x \}$, so that $Q$ has transition function $h_{01}$, we can form the transition function $g_{01} := h_{01} \sigma$ for a bundle $P$, which we will say has trivializations $\varphi_i : P|_{X_i} \to X_i \times G, i = 0,1$.  Then the map
\begin{align*}
s := \psi_0^{-1} \circ \varphi_0 : P|_{X_0} \to Q|_{X_0}
\end{align*}
is a bundle isomorphism, and we obtain a Hecke modification of $Q$ supported at $x$.  It is then not difficult to show that $(P, s)$ is independent of the choices made.

Conversely, given a modification $s : P|_{X_0} \to Q|_{X_0}$, a choice of trivialization $\varphi_1 : P|_{X_1} \to X_1 \times G$ gives a meromorphic section $\tilde{\varsigma}$ of $Q$ over $X_1$ given by
\begin{align*}
y \mapsto s \circ \varphi_1^{-1}(y,e).
\end{align*}
A different choice of $\varphi_1$ amounts to multiplying this section on the right by a holomorphic $G$-valued function on $X_1$, so we get a well-defined class $\varsigma := [ \tilde{\varsigma} ]$ in $\Gr_Q(x)$.  It is clear that these constructions are inverse to each other.
\end{proof}

\subsubsection{Construction of Spaces of Hecke Modifications} \label{SpacesofHms}

Since we have just shown that $\Gr_Q(x)$ is precisely the space of (meromorphic) Hecke modifications of $Q$ supported at $x$, it follows that we have a stratification
\begin{align*}
\Gr_Q(x) = \coprod_{\lambda^\vee \in Y(T)_+} \Gr_Q^{\lambda^\vee}(x).
\end{align*}
As in the proof of Proposition \ref{Hmspace}, a choice of trivialization $\psi_1$ of $Q$ in a neighbourhood of $x$ essentially gives an identification of $\Gr_Q(x)$ with $\Gr_G(x)$, and we pull the stratification back through this identification; as before, this will be independent of the choice of $\psi_1$.

The union
\begin{align*}
\Gr_Q^{\lambda^\vee} := \coprod_{x \in X} \Gr_Q^{\lambda^\vee}(x),
\end{align*}
is thus the set of all Hecke modifications of $Q$ of type $\lambda^\vee$.  It can be given the structure of a fibre bundle over $X$ with fibre isomorphic to $\Gr_G^{\lambda^\vee}$ as follows.  There is an obvious projection map $\pi : \Gr_Q^{\lambda^\vee} \to X$ whose fibre over $x \in X$ is
\begin{align*}
\left(\Gr_Q^{\lambda^\vee} \right)_x = \pi^{-1}(x) := \Gr_Q^{\lambda^\vee}(x).
\end{align*}
Suppose $U \subseteq X$ is an open set over which we have a coordinate $z : U \to z(U) \subseteq \C$ and a trivialization $\psi : Q|_U \arsim U \times G$.  We obtain a bijection $\coprod_{x \in U} \Gr_Q^{, \lambda^\vee}(x) \arsim U \times \Gr_G^{\lambda^\vee}$ as follows.  If $\varsigma \in \Gr_Q^{\lambda^\vee}(x)$ and $\tilde{\varsigma}$ is a locally defined meromorphic section of $Q$ representing $\varsigma$, then we map
\begin{align*}
\varsigma \mapsto \left( x, [p_G \circ \psi \circ \varsigma \circ z^{-1}] \big( z + z(x) \big) \right);
\end{align*}
here $p_G : U \times G \to G$ is the projection map, so that $p_G \circ \psi \circ \varsigma \circ z^{-1}$ is a meromorphic $G$-valued function defined in a neighbourhood of $0 \in \C$, i.e.\ an element of $LG$, and by $[p_G \circ \psi \circ \varsigma \circ z^{-1}]$, we mean its class in $\Gr_G$, which obviously lies in $\Gr_G^{\lambda^\vee}$.  The inverse is given by
\begin{align*}
(x, [\sigma]) \mapsto \left\{ y \mapsto \psi^{-1} \big( x, \sigma( z(y) - z(x)) \big) \right\},
\end{align*}
where $\sigma \in LG$ is a representative for $[\sigma]$.  It is easy to see that these maps are independent of the choices of representatives.

Suppose now that $V \subseteq X$ is an open set on which we have a coordinate $t : V \to \C$ and a trivialization $\varphi : Q|_V \to V \times G$, then if $(\psi, z), (\varphi, t)$ denote the respectively trivializations of $\Gr_Q^{\lambda^\vee}$, we have for $x \in U \cap V$
\begin{align*}
(\psi, z) \circ (\varphi, t)^{-1}( x, \sigma) = \left( x, g_{UV}\big( z^{-1}(z + z(x)) \big) \sigma \big( t \circ z^{-1}(z + z(x)) - t(x) \big) \right),
\end{align*}
where $g_{UV} : U \cap V \to G$ is the transition function for the trivializations $\psi, \varphi$ of $Q$.  We want to see that this gives a holomorphic map $U \cap V \to \Aut \, \Gr_G^{\lambda^\vee}$.  By Corollary \ref{coordchange}, changes of coordinate act holomorphically, so we may assume that $t = z$.  But since $g_{UV}$ is holomorphic, so is
\begin{align*}
x \mapsto g_{UV}\bigg( z^{-1}\big(z + z(x) \big) \bigg).
\end{align*}
Therefore $\Gr_Q^{\lambda^\vee}$ is indeed a holomorphic fibre bundle over $X$.  We will call $\Gr_Q^{\lambda^\vee}$ the \emph{space of Hecke modifications of $Q$ of type $\lambda^\vee$}.  We see that
\begin{align} \label{dimYQlambda}
\dim \Gr_Q^{\lambda^\vee} = \dim \Gr_G^{\lambda^\vee} + 1 = 2 \langle \lambda^\vee, \rho \rangle + 1 = \ell\big( t(\lambda^\vee) \big) + 1.
\end{align}

\subsection{Construction of a Universal Hecke Modification of a Fixed Type}

Fix a $G$-bundle $Q$ over $X$ and a dominant cocharacter $\lambda^\vee \in Y(T)_+$.  In this subsection, we give a construction of a universal Hecke modification $\QQ(\lambda^\vee)$ of $Q$ of type $\lambda^\vee$.  This will be a $G$-bundle over $X \times \Gr_Q^{\lambda^\vee}$ of which we will demand two properties which we now explain.

Let $p : X \times \Gr_Q^{\lambda^\vee} \to X$ and $q : X \times \Gr_Q^{\lambda^\vee} \to \Gr_Q^{\lambda^\vee}$ denote the respective projections.  We also have the projection $\pi : \Gr_Q^{\lambda^\vee} \to X$ from Section \ref{SpacesofHms}.  Therefore, there is a map $\1_X \times \pi : X \times \Gr_Q^{\lambda^\vee} \to X \times X$, and we will define the closed subset $\Gamma \subseteq X \times \Gr_Q^{\lambda^\vee}$ as the preimage of the diagonal $\Delta \subseteq X \times X$:
\begin{align*}
\Gamma := (\1_X \times \pi)^{-1}(\Delta).
\end{align*}
Since $\Delta$ is a divisor on $X \times X$, $\Gamma$ is a divisor on $X \times \Gr_Q^{\lambda^\vee}$.  Let us denote its complement by
\begin{align*}
X_0^\Gr = X_0^{\Gr_Q^{\lambda^\vee}} := X \times \Gr_Q^{\lambda^\vee} \setminus \Gamma.
\end{align*}
The first property we require of $\QQ(\lambda^\vee)$ is for there to be an isomorphism
\begin{align} \label{Qpullback}
\mu : \QQ(\lambda^\vee)|_{X_0^\Gr} \arsim p^*Q|_{X_0^\Gr}.
\end{align}

If $\varsigma \in \Gr_Q^{\lambda^\vee}$, then we will denote by $(Q^\varsigma, s_\varsigma)$ the Hecke modification of $Q$ corresponding to $\varsigma \in \Gr_Q^{\lambda^\vee}$.  The second property we will want $\QQ(\lambda^\vee)$ to satisfy is that of a universal family of modifications in the sense that if $\QQ(\lambda^\vee)_\varsigma := \QQ(\lambda^\vee)|_{X \times \varsigma}$ and $\mu_\varsigma := \mu|_{X \times \varsigma}$, then
\begin{align} \label{univfamily}
\big( \QQ(\lambda^\vee)_\varsigma, \mu_\varsigma \big) \cong (Q^\varsigma, s_\varsigma)
\end{align}
for all $\varsigma \in \Gr_Q^{\lambda^\vee}$.  A bundle $\QQ(\lambda^\vee)$ satisfying (\ref{Qpullback}) and (\ref{univfamily}) will be called a \emph{universal Hecke modification of $Q$ of type $\lambda^\vee$}.

First, we prove a local uniqueness property.

\begin{lem} \label{localuniversalbundle}
Suppose $B \subseteq X \times \Gr_Q^{\lambda^\vee}$ is an open set over which there exist $G$-bundles $\QQ_1$ and $\QQ_2$ and isomorphisms $\mu_1$ and $\mu_2$ as in (\ref{Qpullback}), defined over $B \cap X_0^\Gr = B \setminus \Gamma$, satisfying (\ref{univfamily}) for all $\varsigma \in q(B)$.  Then there exists a unique isomorphism $a : \QQ_1 \to \QQ_2$ such that
\begin{align}  \label{localdiagram}
\xymatrixrowsep{3pc}\xymatrixcolsep{1pc}
\vcenter{\vbox{
    \xymatrix{
    \QQ_1 |_{B \setminus \Gamma} \ar[rr]^a \ar[dr]_{\mu_1} & & \QQ_2 |_{B \setminus \Gamma} \ar[dl]^{\mu_2} \\
    & p^*Q |_{B \setminus \Gamma} & }
    }}
\end{align}
commutes.
\end{lem}

\begin{proof}
We will note that (\ref{localdiagram}) determines $a$ on an open dense set, so by continuity any such isomorphism will necessarily be unique.  The problem is thus reduced to defining $a$.

If $B \subseteq X_0^\Gr$, then there is virtually nothing to prove.  Let $D$ be a neighbourhood of a point in $B \cap \Gamma$.  By taking intersections if necessary, we may take $D$ to be such that it is contained in a set of the form $U \times \pi^{-1}(U)$, where $U \subseteq X$ is an open set over which there exists a trivialization $\psi : Q|_U \arsim U \times G$.  This $\psi$ will induce a trivialization $\overline{\psi} : p^*Q |_{U \times \Gr_Q^{\lambda^\vee}} \arsim U \times \Gr_Q^{\lambda^\vee} \times G$.  We may assume that $D$ is small enough so that we may choose trivializations $b_i : \QQ_i|_D \to D \times G, i = 1,2$.  Then we may write
\begin{align*}
\overline{\psi} \circ \mu_i \circ b_i^{-1} = \1 \times L_{\eta_i}
\end{align*}
for some holomorphic $G$-valued functions $\eta_i : D \setminus \Gamma \to G$.  By (\ref{univfamily}), if we fix $\varsigma \in q(D)$, then if $D_\varsigma := p^{-1}(x) \cap B$, the function $\eta_i^\varsigma : D_\varsigma \setminus \{ \pi(\varsigma) \} \to G$ gives a representative for the twist in the transition function which yields the Hecke modification $\varsigma$.  Since we are assuming that both $\QQ_1$ and $\QQ_2$ satisfy (\ref{univfamily}), it follows that $(\eta_2^\varsigma)^{-1} \eta_1^\varsigma$ gives a holomorphic function on all of $D_\varsigma$.  This implies that $\eta_2^{-1} \eta_1$ is bounded near $D \cap \Gamma$ and so extends to a holomorphic function $D \to G$.  We may thus define an isomorphism $a_D : \QQ_1|_D \to \QQ_2|_D$ by
\begin{align} \label{adef}
a_D := b_2^{-1} \circ (\1_D \times L_{\eta_2^{-1} \eta_1}) \circ b_1.
\end{align}
It is easy to check that we get a commutative diagram (\ref{localdiagram}) over $D \setminus \Gamma$.

Verification that $a_D$ is well-defined is also relatively straightforward.  We made choices of trivializations $\psi$ and $b_1, b_2$.  A change of the trivialization $\psi$ replaces $\eta_i$ by $\nu \eta_i$ for some holomorphic $G$-valued function $\nu : D \to G$.  Then $(\nu \eta_2)^{-1} (\nu \eta_1) = \eta_2^{-1} \eta_1$ and this does not affect the definition of $a_D$.  A change in the trivialization $b_i$ means it is replaced by $(\1 \times \tau_i) \circ b_i$ for some holomorphic $\tau_i : D \to G, i = 1,2$.  Then $\eta_i$ is replaced by $\eta_i \tau_i^{-1}$ and so making these replacements yields the same expression in (\ref{adef}).  Finally, we can cover $B \cap \Gamma$ by such sets $D$ and we get $a : \QQ_1 \to \QQ_2$ defined on all of $B$.
\end{proof}

\begin{thm} \label{univmodification}
There exists a universal Hecke modification $\QQ(\lambda^\vee)$ of $Q$ of type $\lambda^\vee$ for any $G$-bundle $Q$ over $X$ and dominant cocharacter $\lambda^\vee \in Y(T)_+$.
\end{thm}

\begin{proof}
In view of Lemma \ref{localuniversalbundle}, it suffices to construct bundles $\QQ_\alpha$ over open subsets $X_\alpha^\Gr \subseteq X \times \Gr_Q^{\lambda^\vee}$ satisfying (\ref{Qpullback}) and (\ref{univfamily}) for an open cover $\{ X_\alpha^\Gr \}$ of $X \times \Gr_Q^{\lambda^\vee}$.  For in this case, we will obtain isomorphisms $a_{\alpha \beta} : \QQ_\alpha|_{X_{\alpha \beta}^\Gr} \arsim \QQ_\beta|_{X_{\alpha \beta}^\Gr}$, where $X_{\alpha \beta}^\Gr := X_\alpha^\Gr \cap X_\beta^\Gr$.  The cocycle condition on the $a_{\alpha \beta}$ will be satisfied by the uniqueness statement of the Lemma.  Thus, we obtain a bundle on $X \times \Gr_Q^{\lambda^\vee}$ with the required properties.

First, we observe that taking $p^*Q$ over $X_0^\Gr$ gives the required bundle on $X_0^\Gr$.  It remains to show that $\Gamma$ can be covered by open sets over which we have local universal bundles.  Let $\{ U_\alpha \}$ be an open cover of $X$ so that over each $U_\alpha$ we have a coordinate $z_\alpha : U_\alpha \to z_\alpha(U_\alpha) \subseteq \C$ and a trivialization $\psi_\alpha : Q|_{U_\alpha} \arsim U_\alpha \times G$.  Set
\begin{align*}
X_\alpha^\Gr := U_\alpha \times \pi^{-1}(U_\alpha).
\end{align*}
The $X_\alpha^\Gr$ cover $\Gamma$ and hence $\{ X_\alpha^\Gr \} \cup \{ X_0^\Gr \}$ gives an open cover of $X \times \Gr_Q^{\lambda^\vee}$.  We will construct $\QQ(\lambda^\vee)$ by giving a bundle over each open set in a refinement of this cover.

There are isomorphisms $(z_\alpha, \psi_\alpha) : \pi^{-1}(U_\alpha) \arsim U_\alpha \times \Gr_Q^{\lambda^\vee}$ as in Section \ref{SpacesofHms}, so we obtain
\begin{align*}
m _\alpha := \1 \times (z_\alpha, \psi_\alpha) : X_\alpha^\Gr \arsim U_\alpha \times U_\alpha \times \Gr_G^{\lambda^\vee}.
\end{align*}
It is clear that under this isomorphism $\Gamma \cap X_\alpha^\Gr$ corresponds to $\Delta_{U_\alpha} \times \Gr_G^{\lambda^\vee}$.

Let $V_w \subseteq \Gr_G^{\lambda^\vee}$ be an open set as in Theorem \ref{Bruhatdecomp}(b); then if we set
\begin{align*}
X_{\alpha, w}^\Gr := m_\alpha^{-1}( U_\alpha \times U_\alpha \times V_w ),
\end{align*}
it follows that
\begin{align*}
X_\alpha^\Gr = \bigcup_{w \in \lfloor W/W_{t(\lambda^\vee)} \rfloor} X_{\alpha, w}^\Gr.
\end{align*}
We will construct a bundle $\QQ_{\alpha, w}$ on each $X_{\alpha, w}^\Gr$ with the required properties.

Corollary \ref{LGsection} provides us with a section
\begin{align*}
f_w : V_w \to LG
\end{align*}
which allows us to define a function $g_w : (U_\alpha \times U_\alpha \setminus \Delta_{U_\alpha}) \times V_w \to G$ by
\begin{align*}
g_w(x, y, \sigma) = f_w(\sigma)\big(z_\alpha(x) - z_\alpha(y) \big).
\end{align*}
The same Corollary ensures that $g_w$ is well-defined.  We now define a bundle $\QQ_{\alpha, w}$ on $X_{\alpha, w}^\Gr$ by taking trivializations
\begin{align*}
& \varphi_{\alpha,w,0} : \QQ_{\alpha, w}|_{X_{\alpha, w}^\Gr \setminus \Gamma } \to X_{\alpha, w}^\Gr \setminus \Gamma \times G, & & \varphi_{\alpha,w,1} : \QQ_{\alpha, w} \to X_{\alpha, w}^\Gr \times G,
\end{align*}
and with transition function $g_w \circ m_\alpha$, i.e.
\begin{align*}
\varphi_{\alpha,w,0} \circ \varphi_{\alpha,w,1}^{-1} (x, \varsigma,g) = \big(x, \varsigma, g_w \circ m_\alpha(x,\varsigma) g \big).
\end{align*}
Observe that if we think of a fixed $\varsigma \in \pi^{-1}(U_\alpha)$ as a meromorphic section of $Q$, then
\begin{align*}
g_w \circ m_\alpha(x,\varsigma) = \psi_\alpha \circ \tilde{\varsigma}(x),
\end{align*}
where the choice of representative $\tilde{\varsigma}$ section is determined by the section $f_w$ over $V_w$.  One will notice that $\QQ_{\alpha, w}$ is a trivial bundle, but what is important is the relationship to $p^*Q$.  We now show that $\QQ_{\alpha,w}$ satisfies (\ref{Qpullback}) and (\ref{univfamily}).

The trivialization $\psi_\alpha$ yields one $\overline{\psi}_\alpha : p^*Q|_{X_\alpha^\Gr} \arsim X_\alpha^\Gr \to X_\alpha^\Gr \times G$, which we may restrict to $X_{\alpha,w}^\Gr \setminus \Gamma$, and so we can define $\mu_{\alpha, w} : \QQ_{\alpha,w}|_{X_{\alpha, w}^\Gr \setminus \Gamma} \arsim p^*Q|_{X_{\alpha, w}^\Gr \setminus \Gamma}$ by the composition
\begin{align*}
\mu_{\alpha, w} := \overline{\psi}_\alpha^{-1} \circ \varphi_{\alpha,w,0}.
\end{align*}
This is (\ref{Qpullback}).  Now, fix $\varsigma \in \pi^{-1}(U_\alpha)$, say with $x := \pi(\varsigma)$, and consider $\QQ_{\alpha,w}|_{U_\alpha \times \varsigma}$.  Take the trivialization $\varphi_{\alpha,w,1}|_{U_\alpha \times \varsigma} : \QQ_{\alpha,w}|_{U_\alpha \times \varsigma} \arsim U_\alpha \times G$.  We have an isomorphism
\begin{align*}
\mu_{\alpha,w}|_{U_\alpha \times \varsigma} : \QQ_{\alpha,w}|_{U_\alpha \setminus \{ x \} \times \varsigma} \arsim Q|_{U_\alpha \setminus \{ x \} }
\end{align*}
and the meromorphic section
\begin{align*}
y \mapsto \mu_{\alpha,w}|_{U_\alpha \times \varsigma} \circ \varphi_{\alpha,w,1}|_{U_\alpha \times \varsigma}^{-1} (y,e) = \psi_\alpha^{-1} \big( y, g_w \circ m_\alpha(y, \varsigma) \big)
\end{align*}
yields the class of $\varsigma$ by the remarks in the previous paragraph.  This proves the property (\ref{univfamily}).
\end{proof}

\subsubsection{Multiple Modifications} \label{multiplemodifications}

Suppose $\varsigma \in \Gr_Q^{\lambda^\vee}, \varpi \in \Gr_Q^{\mu^\vee}$ are Hecke modifications of types $\lambda^\vee, \mu^\vee \in Y(T)_+$ supported at the respective distinct points $x$ and $y$.  Consider the modification $(Q^\varsigma, s_\varsigma)$ obtained from $\varsigma$.  Then since $\varpi$ is supported away from $x$, viewing it as a class of a meromorphic section of $Q$, the composition $s_\varsigma^{-1} \circ \varpi$ gives a meromorphic section of $Q^\varsigma$ of the same type, so we may think of $\varpi$ as an element of $\Gr_{Q^\varsigma}^{\mu^\vee}$, and we get a bundle $(Q^\varsigma)^\varpi$ and an isomorphism
\begin{align*}
(Q^\varsigma)^\varpi|_{X_0} \xrightarrow{ s_\varpi \circ s_\varsigma} Q|_{X_0},
\end{align*}
where $X_0 := X \setminus \{ x, y \}$.  Since the transition function for $(Q^\varsigma)^\varpi$ can be taken to be on the disjoint union of punctured discs, by writing down a relation between the transition function for $(Q^\varsigma)^\varpi$ and that of $Q$ as in (\ref{tfreln}), one can see that there is a unique isomorphism $\alpha : (Q^\varsigma)^\varpi \to (Q^\varpi)^\varsigma$ such that the diagram
\begin{align*}
\xymatrixrowsep{3pc}\xymatrixcolsep{1pc}
\xymatrix{
(Q^\varsigma)^\varpi|_{X_0} \ar[rr]^\alpha \ar[dr]_{s_\varsigma \circ s_\varpi} & & (Q^\varpi)^\varsigma |_{X_0} \ar[dl]^{s_\varpi \circ s_\varsigma} \\
& Q|_{X_0} & }
\end{align*}
commutes.  Essentially, we are saying that if Hecke modifications are supported at distinct points, then the order in which they are performed does not matter.

For $k \geq 2$, consider the $k$-fold product $X^k$ and the projection $p_i : X^k \to X$ onto the $i$th factor.  For $i < j$, we then get the projections onto two factors $(p_i, p_j) : X^k \to X \times X$.  We let $\Delta_{ij} := (p_i, p_j)^{-1}(\Delta)$, and set
\begin{align*}
\Delta_k := \bigcup_{1 \leq i < j \leq k} \Delta_{ij} \subseteq X^k = \{ (x_1 \ldots, x_k) \in X^k \, : \, x_i = x_j \text{ for some } 1 \leq i,j \leq k \},
\end{align*}
so that $\Delta_k$ is the ``fat'' diagonal.  Note that $\Delta_k$ is a divisor on $X^k$.

Let $\lambda_1^\vee, \ldots, \lambda_k^\vee \in Y(T)_+$.  If $\pi_i : \Gr_Q^{\lambda_i^\vee} \to X$ is the projection map, then we can form the product map
\begin{align*}
\pi_1 \times \cdots \times \pi_k : \Gr_Q^{\lambda_1^\vee} \times \cdots \times \Gr_Q^{\lambda_k^\vee} \to X^k,
\end{align*}
and so we may consider the divisor
\begin{align*}
\widetilde{\Delta}_k := (\pi_1 \times \cdots \times \pi_k)^{-1}(\Delta_k) \subseteq \Gr_Q^{\lambda_1^\vee} \times \cdots \times \Gr_Q^{\lambda_k^\vee}.
\end{align*}
We will let
\begin{align*}
\Gr_{Q,0}^{\lambda_1^\vee, \cdots, \lambda_k^\vee} := \Gr_Q^{\lambda_1^\vee} \times \cdots \times \Gr_Q^{\lambda_k^\vee} \setminus \widetilde{\Delta}_k,
\end{align*}
so that $\Gr_{Q,0}^{\lambda_1^\vee, \cdots, \lambda_k^\vee}$ consists of $k$-tuples of Hecke modifications of $Q$ supported at distinct points.  Set
\begin{align*}
\Gamma_k := (\1_X \times \pi_1 \times \cdots \times \pi_k)^{-1} ( \Delta_{k+1}) \subseteq X \times \Gr_{Q,0}^{\lambda_1^\vee, \cdots, \lambda_k^\vee}.
\end{align*}
We will let $p : X \times \Gr_{Q,0}^{\lambda_1^\vee, \cdots, \lambda_k^\vee} \to X$ and $q_i : X \times \Gr_{Q,0}^{\lambda_1^\vee, \cdots, \lambda_k^\vee} \to \Gr_Q^{\lambda_i^\vee}$ denote the relevant projection maps.  A \emph{universal sequence of Hecke modifications of $Q$ supported at distinct points} will be a $G$-bundle $\QQ(\lambda_1^\vee, \ldots, \lambda_k^\vee)$ over $X \times \Gr_{Q,0}^{\lambda_1^\vee, \cdots, \lambda_k^\vee}$ and an isomorphism
\begin{align} \label{Qpullbackmultiple}
\mu : \QQ(\lambda_1^\vee, \ldots, \lambda_k^\vee) |_{X_0^\Gr} \arsim p^*Q |_{X_0^\Gr},
\end{align}
where ${X_0^\Gr} := X \times \Gr_{Q,0}^{\lambda_1^\vee, \cdots, \lambda_k^\vee} \setminus \Gamma_k$, for which
\begin{align} \label{univfamilymultiple}
\big( \QQ(\lambda_1^\vee, \ldots, \lambda_k^\vee)_{(\varsigma_1, \ldots, \varsigma_k)}, \mu_{(\varsigma_1, \ldots, \varsigma_k)} \big) \cong ( Q^{\varsigma_1 \cdots \varsigma_k}, s_{\varsigma_k} \circ \cdots \circ s_{\varsigma_1}),
\end{align}
for all $(\varsigma_1, \ldots, \varsigma_k) \in \Gr_{Q,0}^{\lambda_1^\vee, \cdots, \lambda_k^\vee}$.  One has analogous results as above.

\begin{lem} \label{localunivbundlemultiple}
Suppose $B \subseteq X \times \Gr_{Q,0}^{\lambda_1^\vee, \cdots, \lambda_k^\vee}$ is an open set over which there exist $G$-bundles $\QQ_1$ and $\QQ_2$ and isomorphisms $\mu_1$ and $\mu_2$ as in (\ref{Qpullbackmultiple}), defined over $B \setminus \Gamma_k$, satisfying (\ref{univfamilymultiple}) for all $\varsigma \in q(B)$.  Then there exists a unique isomorphism $a : \QQ_1 \to \QQ_2$ such that
\begin{align*}
\xymatrixrowsep{3pc}\xymatrixcolsep{1pc}
\vcenter{\vbox{
    \xymatrix{
    \QQ_1 |_{B \setminus \Gamma_k} \ar[rr]^a \ar[dr]_{\mu_1} & & \QQ_2 |_{B \setminus \Gamma_k} \ar[dl]^{\mu_2} \\
    & p^*Q |_{B \setminus \Gamma_k} & }
    }}
\end{align*}
commutes.
\end{lem}

As before, we need only construct the isomorphism in neighbourhoods of points $(x, \varsigma_1, \ldots, \varsigma_k)$ where $x = \pi_i(\varsigma_i)$ for some $1 \leq i \leq k$; but then $x \neq \pi_j(\varsigma_j)$ for $j \neq i$.  Thus, the same argument can be used as for Lemma \ref{localuniversalbundle}.

\begin{thm}
There exists a universal sequence of Hecke modifications $\QQ(\lambda_1^\vee, \ldots, \lambda_k^\vee)$ at distinct points.
\end{thm}

\begin{proof}
As in the proof of Theorem \ref{univmodification}, the Lemma makes it sufficient to construct bundles with the required properties locally.  Inductively, we may assume that a bundle $\QQ(\lambda_1^\vee, \ldots, \lambda_{k-1}^\vee)$ with the required properties exists.  Observe that the projection, which we will call $q$, from $\Gr_{Q,0}^{\lambda_1^\vee, \ldots, \lambda_k^\vee}$ which omits the last factor lands in $\Gr_{Q,0}^{\lambda_1^\vee, \ldots, \lambda_{k-1}^\vee}$ (since if $k$ Hecke modifications are supported at distinct points, then so are the first $k-1$ of them).  Over the open set
\begin{align*}
X_{0,k}^\Gr := \left\{ (x, \varsigma_1, \ldots, \varsigma_k)  \in X \times \Gr_{Q,0}^{\lambda_1^\vee, \ldots, \lambda_k^\vee} \, : \, x \neq \pi_k(\varsigma_k) \right\}
\end{align*}
we may take the bundle to be $q^*\QQ(\lambda_1^\vee, \ldots, \lambda_{k-1}^\vee)$.  Thus, we need only construct locally defined bundles over points of the form $(x, \varsigma_1, \ldots, \varsigma_k)$, where $x = \pi_k(\varsigma_k)$, but this can be done as in the proof of Theorem \ref{univmodification}.
\end{proof}

\subsubsection{Symmetric Products}

Since the order in which we introduce Hecke modifications does not affect the resulting bundle, if we introduce modifications only of the same type, it makes sense that the effective parameter space is a symmetric product.  We formalize this idea in this subsection.  We will suppose that the modifications in question are all of the same type, i.e.\ $\lambda_1^\vee = \cdots = \lambda_k^\vee = \lambda^\vee$.  Then the symmetric group $\mathfrak{S}_k$ acts freely on $\Gr_{Q,0}^{\lambda^\vee, \cdots, \lambda^\vee}$; we will denote the quotient by
\begin{align*}
\left( \Gr_{Q,0}^{\lambda^\vee} \right)^{(k)} := \mathfrak{S}_k \backslash \Gr_{Q,0}^{\lambda^\vee, \cdots, \lambda^\vee},
\end{align*}
and can think of it as unordered $k$-tuples of Hecke modifications supported at distinct points.  Observe that the projection map $\pi \times \cdots \times \pi : \Gr_{Q,0}^{\lambda^\vee, \cdots, \lambda^\vee} \to X^k$ then induces a map
\begin{align*}
\pi^{(k)} : \left( \Gr_{Q,0}^{\lambda^\vee} \right)^{(k)} \to X^{(k)}
\end{align*}
to the $k$th symmetric product of $X$, which we may identify with the space of effective degree $k$ divisors on $X$.  Since we are considering tuples of modifications with support at distinct points, it follow that the image is the open set of $X^{(k)}$ consisting of reduced divisors.

The $\mathfrak{S}_k$-action extends to one on $X \times \Gr_{Q,0}^{\lambda^\vee, \cdots, \lambda^\vee}$, by taking the trivial action on the first factor.  Fix $\nu \in \mathfrak{S}_k$ and consider the bundle $\nu^* \QQ(\lambda^\vee, \ldots, \lambda^\vee)$.  Then we obtain an isomorphism
\begin{align*}
\nu^* \QQ(\lambda^\vee, \ldots, \lambda^\vee) |_{ X_0^\Gr } \xrightarrow{ \nu^* \mu } \nu^* p^* Q|_{ X_0^\Gr } = (p \circ \nu)^* Q|_{ X_0^\Gr } = p^* Q|_{ X_0^\Gr }.
\end{align*}
Furthermore, since the order in which Hecke modifications at distinct points are introduced is immaterial, $\nu^* \QQ(\lambda^\vee, \ldots, \lambda^\vee)$ also satisfies (\ref{univfamilymultiple}).  Hence Lemma \ref{localunivbundlemultiple} yields a unique isomorphism
\begin{align*}
a_\nu : \QQ(\lambda^\vee, \ldots, \lambda^\vee) \to \nu^* \QQ(\lambda^\vee, \ldots, \lambda^\vee)
\end{align*}
with the appropriate commutation properties, and hence a diagram
\begin{align*}
\xymatrix{
\QQ(\lambda^\vee, \ldots, \lambda^\vee) \ar[r]^{a_\nu} \ar[d] & \nu^* \QQ(\lambda^\vee, \ldots, \lambda^\vee) \ar[r] \ar[d] & \QQ(\lambda^\vee, \ldots, \lambda^\vee) \ar[d] \\
X \times \Gr_{Q,0}^{\lambda^\vee, \ldots, \lambda^\vee} \ar@{=}[r] & X \times \Gr_{Q,0}^{\lambda^\vee, \ldots, \lambda^\vee} \ar[r]_\nu & X \times \Gr_{Q,0}^{\lambda^\vee, \ldots, \lambda^\vee}. }
\end{align*}
This justifies the following statement.

\begin{prop}
There is an $\mathfrak{S}_k$-action on the bundle $\QQ(\lambda^\vee, \ldots, \lambda^\vee)$ over the action on $X \times \Gr_{Q,0}^{\lambda^\vee, \ldots, \lambda^\vee}$.  Thus $\QQ(\lambda^\vee, \ldots, \lambda^\vee)$ descends to a bundle $\QQ(\lambda^\vee)^{(k)}$ over $\left( \Gr_{Q,0}^{\lambda^\vee} \right)^{(k)}$, possessing properties analogous to those of (\ref{Qpullbackmultiple}) and (\ref{univfamilymultiple}).
\end{prop}

Suppose now that we are in the situation of Section \ref{multiplemodifications} with $k = \sum_{i=1}^m k_i$ and
\begin{align*}
\lambda_1^\vee = \cdots = \lambda_{k_1}^\vee, \lambda_{k_1 +1}^\vee = \cdots = \lambda_{k_1 + k_2}^\vee, \ldots, \lambda_{k_1 + \cdots + k_{m-1} + 1}^\vee = \cdots = \lambda_k^\vee.
\end{align*}
Then $\Gr_{Q,0}^{\lambda_1^\vee, \ldots, \lambda_k^\vee}$ admits a free action of $\mathfrak{S}_{k_1} \times \cdots \times \mathfrak{S}_{k_m}$, the quotient of which we will denote by
\begin{align*}
\left( \Gr_{Q,0}^{\lambda_{k_1}^\vee, \ldots, \lambda_{k_m}^\vee} \right)^{(k_1, \ldots, k_m)} := \mathfrak{S}_{k_1} \times \cdots \times \mathfrak{S}_{k_m} \backslash \Gr_{Q,0}^{\lambda_1^\vee, \ldots, \lambda_k^\vee}.
\end{align*}
This carries a projection to $X^{(k_1)} \times \cdots \times X^{(k_m)}$ which carries a further projection to $X^{(k)}$.  Then $( \Gr_{Q,0}^{\lambda_{k_1}^\vee, \ldots, \lambda_{k_m}^\vee} )^{(k_1, \ldots, k_m)}$ is precisely the subset of
\begin{align*}
\left( \Gr_{Q,0}^{\lambda_{k_1}^\vee} \right)^{(k_1)} \times \cdots \times \left( \Gr_{Q,0}^{\lambda_{m_1}^\vee} \right)^{(k_m)}
\end{align*}
that maps to the set of reduced divisors in $X^{(k)}$.

The same arguments as above allow us the following.

\begin{prop} \label{symmetricproductbundle}
There is an action of $\mathfrak{S}_{k_1} \times \cdots \times \mathfrak{S}_{k_m}$ on $\QQ(\lambda_1^\vee, \ldots, \lambda_k^\vee)$ which lifts the action on $\Gr_{Q,0}^{\lambda_1^\vee, \ldots, \lambda_k^\vee}$.  There is a universal Hecke modification (i.e.\ bundle satisfying (\ref{Qpullbackmultiple}) and (\ref{univfamilymultiple})) $\QQ(\lambda_{k_1}^\vee, \ldots, \lambda_{k_m}^\vee)^{(k_1, \ldots, k_m)}$ over $\left( \Gr_{Q,0}^{\lambda_{k_1}^\vee, \ldots, \lambda_{k_m}^\vee} \right)^{(k_1, \ldots, k_m)}$.
\end{prop}

\section{Overview of the Wonderful Compactification} \label{OWC}

\subsection{Construction of the Compactification} \label{wcptfn}

We will now assume that $G$ is semisimple of adjoint type.  We will use the same notation as in Section \ref{rootsandweights}; for our choice of Borel subgroup $B = B^+ \subseteq G$, its opposite Borel is denoted $B^-$, and their respective unipotent radicals by $U^+ \subseteq B$ and $U^- \subseteq B^-$.

The compactification $\overline{G}$ of $G$ is constructed as follows.  Let $V$ be a regular irreducible representation of the universal cover $\tilde{G}$ of $G$ with highest weight $\lambda$ (regularity means that $\langle \lambda, \alpha_i^\vee \rangle > 0$ for $1 \leq i \leq l$) and define a map $\psi : G \hookrightarrow \P (\End \, V)$ by
\begin{align*}
g \mapsto [ \tilde{g} ],
\end{align*}
where $\tilde{g} \in \tilde{G}$ is a lift of $g$ and $[ \ ]$ indicates the class in the projectivization.  Then $\overline{G}$ is defined as the closure of $\psi(G)$ in $\P (\End \, V)$.  We will often identify an element of $G$ with its image in $\overline{G}$, which will usually mean abbreviating $\psi(g)$ to $g$.

There is a natural $(G \times G)$-action on $\overline{G}$ given by
\begin{align*}
(g,h) \cdot [\varphi] = [g \varphi h^{-1}].
\end{align*}
We may realize this as separate $G$-actions, which we will call the left and right $G$-actions, corresponding to the action of the first and second factors of $G \times G$, respectively.  We will denote these actions $G \times \overline{G} \to \overline{G}$ by $L$ and $R$, respectively.

\begin{rmk} \label{rightaction}
While we call $R$ the ``right'' action, $h \in G$ acts by $(e,h)$ which takes $[\varphi]$ to $[\varphi h^{-1}]$, so it is in fact a left action, in the sense that we obtain a homomorphism $G \to \Aut \, \GG$, rather than an anti-homomorphism, but we say ``right'' since we mean right multiplication.
\end{rmk}

\subsubsection{The Open Affine Piece} \label{openaffine}

We choose a basis $v_0, \ldots, v_n$ of weight vectors of $V$, say with $v_k$ of weight $\lambda_k$ in such a way that
\begin{enumerate}
\item $v_0$ is a highest weight vector (hence of weight $\lambda$); and
\item $v_1, \ldots, v_l$ are of weights $\lambda_1 = \lambda - \alpha_1, \ldots, \lambda_l = \lambda - \alpha_l$, respectively.
\end{enumerate}
The remaining weights are of the form $\lambda_k = \lambda - \sum n_{ik} \alpha_i$ for some non-negative integers $n_{ik}$.  Let
\begin{align*}
\P_0 = \{ [\varphi] \in \P (\End \, V) \, | \, v_0^*( \varphi v_0) \neq 0 \}.
\end{align*}
If $\varphi = \sum a_{ij} v_i \otimes v_j^*$, then $\P_0$ consists of precisely those $[\varphi]$ with $a_{00} \neq 0$.  Thus $\P_0$ is a standard open affine subset of $\P (\End \, V)$ using the basis $v_i \otimes v_j^*$ of $\End \, V$.  Let $G_0 := \overline{G} \cap \P_0$; then $G_0$ is an open affine subset of
$\overline{G}$.

\begin{lem} \label{affiso} \cite[Lemma 2.6]{EJ}
$G_0 \cap \psi(G) = \psi(U^- T U^+)$.
\end{lem}

We write $Z := \overline{\psi(T)} \cap \P_0$.  Let $t \in T$ and $\tilde{t} \in \tilde{T}$ be a lift, where $\tilde{T}$ is a maximal torus of $\tilde{G}$ mapping onto $T$.  Then if we write $\varphi_k$ for $v_k \otimes v_k^*$, we have
\begin{align*}
\tilde{t} \circ \varphi_k = \tilde{t} \circ v_k \otimes v_k^* = \lambda_k(\tilde{t}) v_k \otimes v_k^* = \lambda_k(\tilde{t}) \varphi_k,
\end{align*}
so
\begin{align} \label{torusmap}
\psi(t) & = [\tilde{t}] = [ \tilde{t} \circ e] = \left[ \sum_{k=0}^n \lambda_k( \tilde{t} ) \varphi_k \right] = \left[ \varphi_0 + \sum_{i=1}^l \frac{ 1 }{\alpha_i(t)} \varphi_i + \sum_{k>l} \prod \frac{ 1 }{ \alpha_i(t)^{n_{ik}} } \varphi_k \right].
\end{align}
Since any two lifts of $t$ differ by an element of $Z(\tilde{G})$, which is precisely the intersection of the $\ker \, \alpha_i, 1 \leq i \leq l$, the values $\alpha_i(t)$ are independent of the choice of lift, so we are justified in dropping the tildes from the notation.

Since $G$ is adjoint, $X(T) = \Lambda_r$, so $\Delta$ is a basis of characters of $T$, hence we may take $z_i = 1/\alpha_i(t), 1 \leq i \leq l$ as coordinates on $T$, and define a map $F : (\C^\times)^l \to Z$ by
\begin{align} \label{Fiso}
(z_1, \ldots, z_l) \mapsto \left[ \varphi_0 + \sum_{i=1}^l z_i \, \varphi_i + \sum_{k>l} \prod z_i^{n_{ik}}  \varphi_k \right].
\end{align}
Clearly, this extends to an isomorphism $\C^l \xrightarrow{\sim} Z$ which we will also denote by $F$.  To lighten notation, we will define
\begin{align*}
p_k = p_k(z_1, \ldots, z_l) := \begin{cases} z_k & 1 \leq k \leq l \\ \prod_{i=1}^l z_i^{n_{ik}} & k > l, \end{cases}
\end{align*}
so that $F$ may be written somewhat more compactly as
\begin{align*}
(z_1, \ldots, z_l) \mapsto \left[ \varphi_0 + \sum_{k=1}^n p_k \varphi_k \right].
\end{align*}

The following statement gives a parametrization of the open affine piece just described that will be basic to the computations in our deformation theory later.
\begin{thm} \label{Aiso}
\cite[Theorem 2.8]{EJ} The map $A : U^- \times U^+ \times Z \xrightarrow{\sim} G_0$ given by
\begin{align*}
(u_-, u_+, t) \mapsto (u_-, u_+) \cdot [t] = [u_- t u_+^{-1}]
\end{align*}
is an isomorphism.  Hence $G_0 \cong \C^{\dim \, G}$.
\end{thm}

The structure of the $(G \times G)$-orbits of $\GG$ has a ready description, and orbit representatives can be chosen in the closure of the torus.

\begin{thm} \label{GxGaction} \cite[Theorem 2.22]{EJ} The compactification $\overline{G}$ is the union of the $(G \times G)$-translates of $G_0$.  In fact, the $(G \times G)$-orbit structure of $\GG$ can be described as follows.  For a subset $I \subseteq \Delta$, set
\begin{align*}
z_I := F(\epsilon_1, \ldots, \epsilon_l),
\end{align*}
where $\epsilon_i = 1$ if $\alpha_i \not\in I$ and $\epsilon_i = 0$ if $\alpha_i \in I$.  Then
\begin{align*}
\GG = \coprod_{I \subseteq \Delta} (G \times G) \cdot z_I,
\end{align*}
so $\GG$ is the union of $2^l$ $(G \times G)$-orbits.  Also, $G \subseteq \GG$ corresponds to the (open) orbit of $z_\emptyset = e$.
\end{thm}

\subsection{The Infinitesimal Action on $T\GG$}

\subsubsection{The Action at a Point}

We now fix a point $z = F(z_1, \ldots, z_l) \in Z$ and give an explicit description of the differentials $dR, dL : \g \to T_z \GG$.  The isomorphism $A : U^- \times U^+ \times Z \to G_0$ of Theorem \ref{Aiso} yields an isomorphism of tangent spaces
\begin{align*}
dA_{(e,e,z)} : T(U^- \times U^+ \times Z) = \u^- \oplus \u^+ \oplus \C^l \to T_z G_0 = T_z \GG,
\end{align*}
which we may use to obtain a basis for $T_z \GG$.  We will shorten to $R, L$ the maps $R_z, L_z : G \to \GG$ given by
\begin{align*}
R_z(g) = g \cdot z, \quad L_z(g) = z \cdot g.
\end{align*}
We will also abbreviate $dA_{(e,e,z)}$ to $dA$.

For $\alpha \in \Phi^+$, we will let $x_\alpha \in \g_\alpha, y_\alpha \in \g_{-\alpha}, h_\alpha = \alpha^\vee \in \t$ be a standard $\sl(2)$ triple, so that $\{ x_\alpha \, | \, \alpha \in \Phi^+ \}$ is a basis for $\u^+$, $\{ y_\alpha \, | \, \alpha \in \Phi^+ \}$ one for $\u^-$, $\{ h_\alpha \, | \, \alpha \in \Delta \}$ one for $\t$.  A basis for $T_z \GG$ is given by
\begin{align} \label{TGbarbasis}
\{ dA(x_\alpha), dA(y_\alpha) \, | \, \alpha \in \Phi^+ \} \cup \{ dA(e_j) \, | \, 1 \leq j \leq l \}.
\end{align}
Observe that
\begin{align} \label{dA}
dA(y_\alpha) & = dR(y_\alpha), & dA(x_\alpha) & = -dL(x_\alpha).
\end{align}

Since $z$ lies in $G_0 \subseteq \P_0$, we will identify $T_z \GG$ with a subspace of $T_z \P_0$, and since $\P_0$ is the standard open affine piece of $\P \, \End \, V$ with $a_{00} \neq 0$, we may identify it with the vector space spanned by $v_i \otimes v_j^*$ with $(i,j) \neq (0,0)$.

We will now be more explicit about the maps $dL, dR$. If $\xi \in \u^-$ or $\u^+$ we obtain
\begin{align*}
dL(\xi) & = \frac{d}{d\epsilon} \bigg|_{\epsilon = 0} z \exp( \epsilon \xi) = \frac{d}{d\epsilon} \bigg|_{\epsilon = 0} \left[ z \circ \exp( \epsilon \xi) \right] = \frac{d}{d\epsilon} \bigg|_{\epsilon = 0} \left[ z + \epsilon (z \circ \xi) \right] = z \circ \xi.
\end{align*}
Here we are identifying $z \in \GG$ with the endomorphism of $V$ it represents (as a tangent vector with reference to the identifications made above). Observe that the terms of $z \circ \xi$ will be of the form $v_i \otimes v_i^* \circ \xi$ and hence there is no $\epsilon$ term involving $v_0 \otimes v_0^*$, so we are using the appropriate affine coordinates in which we are taking the derivative. Now, if $\xi \in \t$, the right action on $\varphi_k$ is
\begin{align*}
\varphi_k \circ \xi = - v_k \otimes \xi(v_k^*) = \langle \xi, \lambda_k \rangle v_k \otimes v_k^* = \langle \xi, \lambda_k \rangle \varphi_k.
\end{align*}
Therefore the infinitesimal action of $\xi \in \t$ on $z$ looks like
\begin{align*}
dL(\xi) & = \frac{d}{d\epsilon} \bigg|_{\epsilon = 0} \left[ \big( 1 + \epsilon \langle \xi, \lambda \rangle \big) \varphi_0 + \sum_{k=1}^n \big( 1 + \epsilon \langle \xi, \lambda_k \rangle \big) p_k \varphi_k \right] \\
& = \frac{d}{d\epsilon} \bigg|_{\epsilon = 0} \left[ \varphi_0 + \sum_{k=1}^n \big[ 1 + \epsilon \big( \langle \xi, \lambda_k \rangle - \langle \xi, \lambda \rangle \big) \big] p_k \varphi_k \right] \\
& = \sum_{k=1}^n \langle \xi, \lambda_k - \lambda \rangle p_k \varphi_k = - \left( \sum_{i=1}^l \langle \xi, \alpha_i \rangle z_i \varphi_i + \sum_{k>l} \left\langle \xi, \sum_{i=1}^l n_{ik} \alpha_i \right\rangle \prod_{i=1}^l z_i^{n_{ik}} \varphi_k \right)
\end{align*}
Note that we need to normalize the coefficient of $\varphi_0 = v_0 \otimes v_0^*$.  We can repeat these calculations for $dR$ as well and sum up our results in the following.

\begin{lem} \label{infaction}
Explicit descriptions of the infinitesimal action can be given by
\begin{align*}
dL( \xi) & = \begin{cases} z \circ \xi & \xi \in \u^- \oplus \u^+ \\
- \left( \sum_{i=1}^l \langle \xi, \alpha_i \rangle z_i \varphi_i + \sum_{k>l} \left\langle \xi, \sum_{i=1}^l n_{ik} \alpha_i \right\rangle \prod_{i=1}^l z_i^{n_{ik}} \varphi_k \right) & \xi \in \t
\end{cases}
\end{align*}
and
\begin{align*}
dR(\xi) & = \begin{cases} \xi \circ z & \xi \in \u^- \oplus \u^+ \\
- \left( \sum_{i=1}^l \langle \xi, \alpha_i \rangle z_i \varphi_i + \sum_{k>l} \left\langle \xi, \sum_{i=1}^l n_{ik} \alpha_i \right\rangle \prod_{i=1}^l z_i^{n_{ik}} \varphi_k \right) & \xi \in \t.
\end{cases}
\end{align*}
\end{lem}

By (\ref{dA}), this gives expressions for $dA(y_\alpha), dA(x_\alpha)$.  Finally, $dA(e_j)$ can be computed as
\begin{align*}
dA(e_j) & = \frac{d}{d\epsilon} \bigg|_{\epsilon = 0} F(z_1, \ldots, z_j + \epsilon, \ldots, z_l) = \varphi_j + \sum_{k > l} \pd{p_k}{z_j} \varphi_k .
\end{align*}

This allows us to write down the image of the infinitesimal action of $\t$ in terms of the basis (\ref{TGbarbasis}) for $T_z \GG$.  We will let $h_i := h_{\alpha_i} = \alpha_i^\vee$ for $1 \leq i \leq l$, so that the $h_1, \ldots, h_l$ give a basis for $\t$.  By Lemma \ref{infaction}, using the fact that for a weight $\mu$, $\mu(h_\alpha) = \langle \alpha^\vee, \mu \rangle$, we get
\begin{align*}
dL(h_i) = dR(h_i) & = - \sum_{j=1}^l a_{ij} z_j \, dA(e_j),
\end{align*}
where $(a_{ij}) = ( \langle \alpha_i^\vee, \alpha_j \rangle)$ is the Cartan matrix.

We would also like to write $dL(y_\alpha)$ in terms of the basis vectors $dA(x_\alpha), dA(y_\alpha)$ and $dA(e_j)$ for $T_z\GG$.  Assume first that $z \in T \subseteq G$; then since $\alpha = \sum_{i=1}^l \langle \lambda_i^\vee, \alpha \rangle \alpha_i$
\begin{align*}
dL( y_\alpha ) = dR \circ \Ad \, z (y_\alpha) = \alpha(z)^{-1} dR( y_\alpha) = \prod_{i=1}^l z_i^{ \langle \lambda_i^\vee, \alpha \rangle} dA( y_\alpha).
\end{align*}
By continuity, the formula must also hold for any $z \in Z$.  We repeat for $dR( x_\alpha)$ and summarize our findings as follows.

\begin{lem} \label{infaction2}
The infinitesimal actions in terms of the bases for $\g$ and $T_z\GG$ chosen above are given by
\begin{align*}
dL(x_\alpha) & = -dA(x_\alpha), & dL(y_\alpha) & = \prod_{i=1}^l z_i^{ \langle \lambda_i^\vee, \alpha \rangle} dA(y_\alpha), & dL(h_i) & = - \sum_{j=1}^l a_{ij} z_j \, dA(e_j),
\end{align*}
and
\begin{align*}
dR(x_\alpha) & = -\prod_{i=1}^l z_i^{ \langle \lambda_i^\vee, \alpha \rangle} dA(x_\alpha), & dR(y_\alpha) & = dA(y_\alpha), & dR(h_i) & = - \sum_{j=1}^l a_{ij} z_j \, dA(e_j).
\end{align*}
\end{lem}

\subsubsection{Weyl Twists}

Let $z \in Z \subseteq \overline{T}$ be as in the previous subsection.  Recall that the Weyl group $W = N_G(T)/T$ acts on $T$ and hence on $\overline{T}$.  However, if $z \not\in T$ and $\nu \in W$ is not the identity, then $\nu \cdot z \not\in G_0$.  In any case, if $w \in N_G(T)$ is a representative for $\nu$, then $\nu \cdot z$ lies in the open affine set $w G_0 w^{-1}$. Composing the isomorphism $A$ with conjugation by $w$ gives an isomorphism $A^\nu : U^- \times U^+ \times Z \to w G_0 w^{-1}$:
\begin{align*}
(u_-, u_+, z) & \mapsto [ w u_- z u_+^{-1} w^{-1}] = \left[ \big( \Ad \, w (u_-) \big) (\nu \cdot z) \big( \Ad \, w (u_+^{-1}) \big) \right].
\end{align*}

Since $\Ad \, w^{-1}$ gives an automorphism of $\u^- \oplus \u^+$, we can take
\begin{align*}
\{ dA^\nu( \Ad \, w^{-1} x_\alpha), dA^\nu( \Ad \, w^{-1} y_\alpha) \, | \, \alpha \in \Phi^+ \} \cup \{ dA^\nu(e_i) \, | \,  1 \leq i \leq l \}
\end{align*}
as a basis of $T_{\nu z} \GG$.

We may now compute the infinitesimal actions as in Lemma \ref{infaction2} using this basis.  We have for $\xi \in \u^- \oplus \u^+$,
\begin{align*}
dL_{\nu \cdot z} ( \xi) = \frac{d}{d \epsilon} \bigg|_{\epsilon = 0} w z w^{-1} \exp( \epsilon \xi) = \frac{d}{d \epsilon} \bigg|_{\epsilon = 0}  w z \exp( \epsilon \Ad \, w^{-1} \xi ) w^{-1}.
\end{align*}
If $\alpha \in \nu \Phi^+$, then $\nu^{-1} \alpha \in \Phi^+$ and $\Ad \, w^{-1} x_\alpha$ is a positive root vector, so
\begin{align*}
dL_{\nu \cdot z} (x_\alpha) = \frac{d}{d \epsilon} \bigg|_{\epsilon = 0} A^\nu( e, \exp( - \epsilon \Ad \, w^{-1} x_\alpha), z ) = - dA^\nu ( \Ad \, w^{-1} x_\alpha).
\end{align*}
Also,
\begin{align*}
dL_{\nu \cdot z} ( y_\alpha) & = \frac{d}{d \epsilon} \bigg|_{\epsilon = 0} w \exp( \epsilon \Ad zw^{-1} y_\alpha ) zw^{-1} = \frac{d}{d \epsilon} \bigg|_{\epsilon = 0} A^\nu \big( \exp( \epsilon \Ad zw^{-1} y_\alpha ), e, z \big) \nonumber \\
& = dA^\nu ( \Ad \, z \circ \Ad \, w^{-1} y_\alpha ) = (\nu^{-1} \alpha)(z)^{-1} dA^\nu ( \Ad \, w^{-1} y_\alpha ) \nonumber \\
& = \prod_{i=1}^l z_i^{ \langle \lambda_i^\vee, \nu^{-1} \alpha \rangle } dA^\nu ( \Ad \, w^{-1} y_\alpha).
\end{align*}

We record this here.

\begin{lem} \label{infactiontwist}
If $z \in Z, \nu \in W$ and $w \in N_G(T)$ is a chosen representative for $\nu$, then the infinitesimal action of the root vectors on $\nu \cdot z$ is given by
\begin{align*}
dL_{\nu \cdot z} ( x_\alpha) & = - dA^\nu ( \Ad \, w^{-1} x_\alpha), & dL_{\nu \cdot z} (y_\alpha) & = \prod_{i=1}^l z_i^{ \langle \lambda_i^\vee, \nu^{-1} \alpha \rangle } dA^\nu ( \Ad \, w^{-1} y_\alpha).
\end{align*}
\end{lem}

\subsubsection{Transposes}

Returning back to the situation where $z \in Z$, we obtain the following mirror images to the infinitesimal action maps, which will be useful in the sequel.

\begin{lem} \label{infactiontranspose}
With respect to the basis dual to that used above the transposes of the infinitesimal actions, $dL^t, dR^t : T^*_{z_I} \GG \to \g^*$ are given by
\begin{align*}
dL^t \big( dA(x_\alpha)^* \big) & = -x_\alpha^*, & dL^t \big( dA(y_\alpha)^* \big) & = \prod_{i=1}^l z_i^{\langle \lambda_i^\vee, \alpha \rangle} y_\alpha^*, & dL^t \big( dA(e_i)^* \big) & = - z_i \sum_{j=1}^l a_{ji} h_j^*.
\end{align*}
and
\begin{align*}
dR^t( dA(x_\alpha)^* ) & = - \prod_{i=1}^l z_i^{\langle \lambda_i^\vee, \alpha \rangle} x_\alpha^*, & dR^t( dA(y_\alpha)^*) & = y_\alpha^*, & dR^t( dA(e_i)^* ) & = - z_i \sum_{j=1}^l a_{ji} h_j^*.
\end{align*}
\end{lem}

Recall that the Killing form $\kappa$ gives a $\Ad$-invariant non-degenerate symmetric pairing on $\g$ and hence gives an $\Ad$-equivariant isomorphism $\tilde{\kappa} : \g^* \to \g$.  For $\alpha \in \Phi^+$, there will be $c_\alpha \in \C^\times$ such that
\begin{align*}
\tilde{\kappa}(x_\alpha^*) = c_\alpha y_\alpha, \quad \tilde{\kappa}(y_\alpha^*) = c_\alpha x_\alpha,
\end{align*}

It is straightforward to see that in the case $g \in G$, the maps
\begin{align*}
\xymatrix{ T_g^*G \ar[r]^{dL_g^t}_{dR_g^t} & \g^* \ar[r]^{\tilde{\kappa}} & \g \ar[r]^{dL_g}_{dR_g} & T_g G. }
\end{align*}
agree.  For in this case $L_g = R_g \circ \Ad \, g$ and hence
\begin{align*}
dL_g \circ \tilde{\kappa} \circ dL_g^t & = dR_g \circ \Ad \, g \circ \tilde{\kappa} \circ (dR_g \circ \Ad \, g)^t = dR_g \circ \Ad \, g \circ \tilde{\kappa} \circ (\Ad \, g^{-1})^* \circ dR_g^t \\
& = dR_g \circ \tilde{\kappa} \circ dR_g^t.
\end{align*}
In fact, this is true at any point of $\GG$.

\begin{prop} \label{lrtranspose}
If $a \in \GG$, then the maps
\begin{align*}
\xymatrix{ T_a^*\GG \ar[r]^{dL^*}_{dR^*} & \g^* \ar[r]^{\tilde{\kappa}} & \g \ar[r]^{dL}_{dR} & T_a \GG }
\end{align*}
agree.
\end{prop}

\begin{proof}
If $a \in Z$, it is straightforward to verify this using the expressions in Lemmas \ref{infaction2} and \ref{infactiontranspose}.  In general, we may write $a = g z h$ for some $g, h \in G, z \in Z$, use the equalities
\begin{align*}
L_a & = L_g \circ L_z \circ L_h, & R_a & = R_h \circ R_z \circ R_g,
\end{align*}
and the fact that $R, L$ commute.
\end{proof}

\subsection{Extension of the Inversion Map to $\overline{G}$}

We show that if $\iota : G \to G$ is the inversion map taking $g$ to $g^{-1}$, then it extends to an involution $\iota : \overline{G} \to \overline{G}$.  The idea is as follows. We consider an open affine set $G_\omega$ of $\GG$, this time built from the lowest weight vector rather than the highest weight vector.  Then we construct an isomorphism of $\C^l \cong Z \subseteq G_0$ onto a subset $W \subseteq G_\omega$, and then show that inversion in $T$ extends to $Z$.  We can then extend it from $G_0$ to $G_\omega$ and from there to all of $\GG$ by using the $(G \times G)$-action on $\GG$.

Let $W$ denote the Weyl group corresponding to $T$ and let $\omega \in W$ be (the unique element) such that $\omega( \Delta ) = - \Delta$; write
\begin{align*}\omega(\alpha_i) = -\alpha_{\omega i},\end{align*}
where we are using $\omega$ also to denote the permutation of the indices.  Note that $\omega^2 = 1$.

Now, $\omega \lambda$ will be a lowest weight vector for $V$.  By relabelling, we may assume that $\lambda_n = \omega \lambda$. Further, for $1 \leq i \leq l$, we may also assume
\begin{align*}
\lambda_{n-i} = \omega \lambda_i = \omega( \lambda - \alpha_i) = \lambda_n + \alpha_{\omega i}.
\end{align*}
We now set
\begin{align*}
\P_\omega := \{ [\varphi] \in \P (\End \, V) \, | \, v_n^*(\varphi v_n) \neq 0 \}
\end{align*}
and let $G_\omega := \GG \cap \P_\omega$.  Then as for Lemma \ref{affiso}, we can show that
\begin{align*}
G_\omega \cap \psi(G) = \psi(U^+ T U^-)
\end{align*}
and that $\nu : U^+ \times U^- \times W \to G_\omega$ given by
\begin{align*}
(v_+, v_-, s) = [v_+ s v_-^{-1}]
\end{align*}
is an isomorphism.

Repeating the calculation of (\ref{torusmap}) above, if $t \in T$ and $\tilde{t} \in \tilde{T}$ is a lift, then
\begin{align*}
\psi(t) & = \left[ \sum_{k < n-l} \prod_{i=1}^l \alpha_{\omega i} (t)^{m_{ik}} \varphi_k  + \sum_{i=1}^l \alpha_{\omega i}(t) \varphi_{n-i} + \varphi_n \right],
\end{align*}
where if $k < n-l$, then $\lambda_k = \lambda_n + \sum m_{ik} \alpha_{\omega i}$ with the $m_{ik} \geq 0$.

We let $W := \overline{\psi(T)} \cap \P_\omega$.  Then as in (\ref{Fiso}), we can construct an isomorphism $H : \C^l \to W$:
\begin{align*}
(w_1, \ldots, w_l) \mapsto \left[ \sum_{k < n-l} \prod_{i=1}^l w_i^{m_{ik}} \varphi_k + \sum_{i=1}^l w_i \varphi_{n-i} + \varphi_n \right].
\end{align*}

We now define $a : \C^l \to \C^l$ by
\begin{align*}(z_1, \ldots, z_l) \mapsto (z_{\omega 1}, \ldots, z_{\omega l}).\end{align*}
Then $\iota := H \circ a \circ F^{-1}$ gives an isomorphism $Z \to W$ and if $t \in T$, we may verify that
\begin{align*}\iota(t) = H \circ a \circ F^{-1}(t) = t^{-1}\end{align*}
by writing
\begin{align*}
t = \left[ \varphi_0 + \sum_{k=1}^n \prod_{i=1}^l \alpha_i(t^{-1})^{n_{ik}} \varphi_k \right],
\end{align*}
and using the definitions of the $n_{ik}, m_{ik}$.

We can extend this map to $\iota : G_0 \cong U^- \times U^+ \times Z \to G_\omega \cong U^+ \times U^- \times W$ by
\begin{align*}\iota( [u_- t u_+^{-1}] ) = [ u_+ \iota(t) u_-^{-1}].\end{align*}
We observe that if $g \in G \cap X_0$, then $\iota(g) = g^{-1}$.

To extend $\iota$ to all of $\GG$, we will need the following lemma, which can be proved in the same way as Proposition 2.25 of \cite{EJ} by replacing $v_0$ by $v_n$ where appropriate.

\begin{lem}
Let $z_I$ be as in Theorem \ref{GxGaction}.  If we define $w_I := \iota(z_I)$, then $(g,h) \in (G \times G)_{z_I}$ if and only if $(h,g) \in (G \times G)_{w_I}$.
\end{lem}

The extension of $\iota$ to all of $\GG$ now proceeds straightforwardly as follows.  For $a \in \GG$, let $(g,h) \in G \times G$ be such that $(g,h) \cdot z_I = a$ (Theorem \ref{GxGaction}).  We define
\begin{align*}
\iota(a) := (h, g) \cdot w_I = (h,g) \cdot \iota(z_I).
\end{align*}
The Lemma above guarantees that this is well-defined.

Observe that if $a \in G$, then $a = (a, e) \cdot e = (a,e) \cdot z_\emptyset$, so
\begin{align*}
\iota(a) = (e,a) \cdot w_\emptyset = (e,a) \cdot e = a^{-1}.
\end{align*}
So $\iota$ does indeed restrict to the inversion map on $G$.  From this it follows that $\iota^2$ is the identity on an open dense subset of $\GG$, so this holds on all of $\GG$; that is, $\iota$ is an involution.  Further, if $\Lambda_{(g,h)}$ denotes the map $\GG \to \GG$ given by the action of $(g,h)$:
\begin{align*}
a \mapsto (g,h) \cdot a,
\end{align*}
then we have $\iota \circ \Lambda_{(g,h)} = \Lambda_{(h,g)} \circ \iota$ on the open dense set $G \subseteq \GG$ and again this must hold on $\GG$.  We now summarize our results.

\begin{prop} \label{Gbarinvextn}
There exists an involution $\iota : \GG \to \GG$ such that if $a \in G$, then
\begin{align*}\iota(a) = a^{-1}\end{align*}
and
\begin{align*}\iota \big( (g,h) \cdot a \big) = (h,g) \cdot \iota(a)\end{align*}
for all $a \in \GG, (g,h) \in G \times G$.
\end{prop}

\section{Deformation Theory Using the Compactifications} \label{UsingTheCompactifications}

\subsection{Constructions}

Here we explain how the wonderful compactification can be used to compactify the fibres of a principal bundle.  We may then map Hecke modifications of the original bundle into this compactified bundle.  The deformation space we want will be the quotient of vector bundles derived from this constructions.  First, we begin by setting down some notational conventions.

\subsubsection{Conventions and Notations} \label{Gbdlnotation}
We will take $X$ and $G$ as before and fix a $G$-bundle $Q$.  If $\psi : Q|_U \xrightarrow{\sim} U \times G$ is a trivialization of $Q$ over an open $U \subseteq X$, we will denote the corresponding section by $b : U \to Q|_U$, so that
\begin{align*}
b(y) = \psi^{-1}(y, e)
\end{align*}
for $y \in U$.  When subscripts are used on $\psi$, they will likewise be appended to the corresponding $b$.

Since our discussion will centre around Hecke modifications, we will often have occasion to trivialize $Q$ on the complement $X_0$ of a point $x \in X$, say via $\psi_0$.  If $X_1$ is a neighbourhood of $x$ over which $\psi_1 : Q|_{X_1} \to X_1 \times G$ is a trivialization, then the transition function $h_{01} : X_{01} \to G$ satisfies
\begin{align*}
\psi_0 \circ \psi_1^{-1}(y,g) = \big( y, h_{01}(y) g \big).
\end{align*}
If we set
\begin{align*}
h_i := p_G \circ \psi_i : Q|_{X_i} \to G,
\end{align*}
where $p_G : X_i \times G \to G$ is the projection, then we have
\begin{align} \label{bdlid}
b_0(y) h_{01}(y) = b_1(y), \quad q = b_i \big( \pi(q) \big) h_i(q), \quad h_i \circ b_j(y) = h_{ij}(y)
\end{align}
for $y \in X, q \in Q$, wherever the expressions make sense.

For a bundle $P$, we will typically call its trivializations $\varphi_i : P|_{X_i} \to X_i \times G$, the corresponding sections $a_i : X_i \to P|_{X_i}$, its transition function $g_{01}$, and the $G$-components $g_i := p_G \circ \varphi_i$.  The relations (\ref{bdlid}) hold mutatis mutandis.

\subsubsection{Compactifying Bundles}

Since the principal bundle $Q$ comes with fibres isomorphic to $G$, we may wish to use the compactification $\GG$ to compactify the fibres of $Q$.  Recall that $\GG$ comes with both a left and right $G$-action (see Section \ref{wcptfn}).  We form the associated bundle $\overline{Q} := Q \times_G \overline{G}$ using the left action, i.e.\ we take
\begin{align*}
(q,a) \sim (q \cdot g, g^{-1} a) = (q \cdot g, L_{g^{-1}} a) .
\end{align*}
The equivalence class of $(q,a)$ will be denoted $[q : a]$.  Then $\overline{Q}$ still admits a right action
\begin{align*}
[q : a] \cdot g = [q : ag] = [q : R_{g^{-1}} a].
\end{align*}
This achieves the compactification of the fibres of $Q$ suggested above.  We have a natural inclusion $Q \hookrightarrow \overline{Q}$ explicitly given by
\begin{align*}
q \mapsto [q : e],
\end{align*}
which is equivariant with respect to the right action.

\begin{rmk}
One will observe that the action of $g \in G$ on $\overline{Q}$ is via $R_{g^{-1}}$.  Recall from Remark \ref{rightaction} that $R$ is really a left action---in the sense that we obtain a homomorphism $G
\to \Aut(\GG)$, and hence one $G \to \Aut(\overline{Q})$, rather than an anti-homomorphism---so putting in the inverse legitimately turns it into a right action which coincides with right multiplication on $G$ and hence with the (right) action of $G$ on $Q$, considered as a subvariety of $\overline{Q}$.
\end{rmk}

By an automorphism of $\overline{Q}$ we mean an equivariant bundle automorphism (i.e., covering the identity) over $X$.  An automorphism of $Q$ determines one of any associated bundle and hence one of $\overline{Q}$.  On the other hand, any equivariant automorphism of $\overline{Q}$ maps $Q$ to $Q$, so restricting to $Q$ allows us to recover the automorphism of $Q$ from which the one of $\overline{Q}$ arises.

\begin{lem} \label{AutQbar}
We may identify
\begin{align*}
\Aut \, Q = \Aut \, \overline{Q}.
\end{align*}
\end{lem}

\subsubsection{Compactification and Hecke Modifications}

Let $X, G, Q, \overline{Q}, x \in X, X_0$ be as above.  Let $\varsigma \in \Gr_Q^{\lambda^\vee}(x) \subseteq \Gr_Q^{\lambda^\vee}$ be a Hecke modification of $Q$ of type $\lambda^\vee$ supported at $x$.  Let $P = P^\varsigma$ be the bundle we obtain and
\begin{align*}s : P|_{X_0} \xrightarrow{\sim} Q|_{X_0}\end{align*}
the given isomorphism (\ref{Hmdefn}).  Recall that a choice of trivialization $\psi_1$ of $Q$ in a neighbourhood $X_1$ of $x$ over which a representative of $\varsigma$ is defined identifies $\varsigma$ with an element of $\Gr_G(x, \lambda^\vee)$.  A choice of trivialization $\varphi_1$ of $P|_{X_1}$ then allows us to write
\begin{align*}\psi_1 \circ s \circ \varphi_1^{-1} = \1_{X_{01}} \times L_\sigma\end{align*}
for some holomorphic $\sigma : X_{01} \to G$ which will then be a representative for $\varsigma$.

Since $X_{01} := X_1 \setminus \{ x \}$ and $\GG$ is complete, $\sigma$ extends uniquely to a holomorphic map $X_1 \to \GG$ which we will also denote by $\sigma$.  This extension allows us to define a map $P \to \overline{Q}$, which we also denote by $s$, by
\begin{align*} 
s(p) = \begin{cases}
\psi_0^{-1} \circ \varphi_0(p) & \text{ if } \pi_P(p) \in X_0 \\
\psi_1^{-1} \circ (\1 \times L_\sigma) \circ \varphi_1(p) &
\text{ if } \pi_P(p) \in X_1, \end{cases}
\end{align*}
where $L_\sigma$ means left multiplication by $\sigma$.  Using (\ref{tfreln}), it is elementary to see that $s$ is well-defined.  Further, we may observe that $s$ extends the isomorphism of $P|_{X_0}$ onto $Q|_{X_0} \subseteq \overline{Q}|_{X_0}$ and is equivariant with respect to the right action, since the trivializations and $\1 \times L_\sigma$ all commute with the right action.

\subsubsection{The Dual Construction} \label{dualconstruction}

Since everything we have done so far is symmetric in $P$ and $Q$, we may turn it all aroundand consider the map $\tau := \sigma^{-1} = \iota \sigma: X_{01} \to G$, which must extend uniquely to a map $\tau : X_1 \to \overline{G}$.  Since the relation between transition functions of $P$ and $Q$ above can be rewritten as $h_{01} = g_{01} \sigma^{-1} = g_{01} \tau$, we may likewise construct a map $t = \iota s : Q \to \overline{P}$, which will be an isomorphism of $Q|_{X_0}$ onto $P|_{X_0}$.

\subsection{Deformation Theory}

\subsubsection{A Local Description of the Deformation Space}

We now consider the deformation theory of the construction given in the preceding subsection.  We begin concretely with a local description of infinitesimal deformations.  Since the map $s : P \to \overline{Q}$ is given locally by $\sigma : X_1 \to \GG$, an infinitesimal deformation of $s$ is determined by one of $\sigma$.  The latter is given by a section of $\sigma^* T\GG$ over $X_1$.  Since what is important is the class of $\sigma$ in $\Gr_G(x, \lambda^\vee)$, trivial changes to $\sigma$ are given by right multiplication by elements of $\G(X_1)$, i.e.\ holomorphic $G$-valued functions over $X_1$.  Therefore, the trivial deformations are precisely those arising from the infinitesimal right action of $\g$ on $\GG$.

Here and later, we will use the following notation:  if $Y$ is a space with structure sheaf $\OO_Y$, then for any (finite-dimensional) vector space $W$, we will write $\OO_Y^W$ for the sheaf of $W$-valued functions on $Y$.

Since $G$ is open and dense in $\GG$ and is acted upon freely by the right action, the infinitesimal right action gives an inclusion of sheaves $dR : \OO^\g_{\GG} \to T\GG$.  Pulling back via $\sigma : X_1 \to \GG$ we get another inclusion of sheaves on $X_1$:
\begin{align} \label{localdefmnmap}
\OO^\g_{X_1} = \sigma^* \OO^\g_{\GG} \xrightarrow{\sigma^* dR} \sigma^* T\GG.
\end{align}

If we denote by $u_\xi$ the fundamental vector field on $\GG$ generated by the right infinitesimal action of $\xi \in \g$, then the image of $(y, \xi) \in X_1 \times \g$ is
\begin{align} \label{inftmlrightaction}
u_{\xi} \big( \sigma(y) \big) = dL_{\sigma(y)} ( \xi).
\end{align}
Therefore, the non-trivial deformations are given by sections of $\Gamma(X_1, \sigma^* T\overline{G})$ modulo those of the form just given; so we may identify our deformation space with
\begin{align} \label{localdefmn}
\Gamma(X_1, \sigma^* T\overline{G})/ \Gamma(X_1, \OO^\g_{X_1}).
\end{align}
Since $\sigma(X_{01}) \subseteq G$, $\sigma^* dR$ is an isomorphism on $X_{01}$ and so the quotient will be supported at $x$.  Since $\sigma$ is holomorphic at $x$, the space of sections will be finite-dimensional.

\subsubsection{Invariant Construction of the Deformation Space}

We now describe the infinitesimal deformation space of in terms of the bundles $P, Q$ and various associated cohomology groups.  In what follows, $v_\xi, w_\xi$ will denote the fundamental vector fields determined by $\xi \in \g$ on $P$ (or $\overline{P}$), and $Q$ (or $\overline{Q}$), respectively.

We now realize the space of infinitesimal deformations of $s$ as the global sections of a torsion sheaf on $X$, depending on $P$ and $Q$, which requires the following construction to define.  If $\pi_P, \pi_{\overline{Q}}$ are the projection maps to $X$, we let
\begin{align*}
VP := \ker d\pi_P, \quad V\overline{Q} := \ker d\pi_{\overline{Q}}
\end{align*}
be the respective vertical tangent bundles.  Then since $\pi_P = \pi_{\overline{Q}} \circ s$, we have $ds(VP) \subseteq V\overline{Q}$ and there is a map of sheaves $ds : VP \to s^* V \overline{Q}$ on $P$.  Explicitly, this map is
\begin{align*}
v_\xi(p) \mapsto \big( p, ds(v_\xi(p) ) \big) = \big( p, w_\xi( s(p) ) \big),
\end{align*}
where we are realizing $s^* V \overline{Q}$ as
\begin{align*}
s^* V \overline{Q} = \{ ( p, w ) \in P \times V\overline{Q} \, | \, w \in V_{s(p)} \overline{Q} \}.
\end{align*}
Since $\pi_P, \pi_{\overline{Q}}$ are $G$-invariant maps, $VP$ and $V\overline{Q}$ are $G$-invariant sub-bundles of $TP$ and $T\overline{Q}$, respectively, so they admit $G$-actions, linear on the fibres, for which $ds$ is equivariant (since $s$ is).  Quotienting, we get a map of sheaves on $X$,
\begin{align*}
VP/G = \ad \, P \xrightarrow{ds/G} s^*V\overline{Q}/G = : E_{PQ}.
\end{align*}
Since $s^* V\overline{Q} = \{ (p,w) \, | \, w \in V_{s(p)} \overline{Q} \}$, we may concretely write $ds/G : \ad \, P \to E_{PQ}$ as
\begin{align*}
[p : \xi] \mapsto \big[ p : w_\xi \big( s(p) \big) \big].
\end{align*}

\subsubsection{Trivializations for $E_{PQ}$} \label{trivializationsEPQ}

We now describe trivializations for the bundle $E_{PQ}$ just constructed to specify the link between the local description of the infinitesimal deformation space and the bundle $E_{PQ}$ in the previous subsection, and also to facilitate later calculations.  We begin by giving trivializations for the bundle $s^* V \overline{Q}$ on $P$, and then using these to obtain trivializations for $E_{PQ}$.  We will make heavy use of the notation of Section \ref{Gbdlnotation} and will often write $E$ for $E_{PQ}$ when no confusion will arise.

Since $s|_{X_0}$ is an isomorphism of $P|_{X_0}$ onto the open subset $Q|_{X_0}$ of $\overline{Q}|_{X_0}$, we get an isomorphism $ds|_{X_0} : VP|_{X_0} \xrightarrow{ \sim} s^* V\overline{Q}|_{X_0}$, and so the trivialization $P|_{X_0} \times \g \to VP|_{X_0}$ extends to one $\alpha_0^{-1} : P|_{X_0} \times \g \to s^* V\overline{Q}|_{X_0}$ by
composing with $ds$:
\begin{align*}
(p, \xi) \mapsto \big( p, w_\xi(s(p)) \big).
\end{align*}

The vertical tangent bundle on the trivial $\GG$-bundle $X_1 \times \GG$ is given by
\begin{align*}
V(X_1 \times \GG) = p_{\GG}^* T\GG,
\end{align*}
where $p_{\GG} : X_1 \times \GG \to \GG$ is the projection.  Then the trivialization $\psi_1 : \overline{Q}|_{X_1} \to X_1 \times \GG$ induces an isomorphism $V\overline{Q}|_{X_1} \xrightarrow{\sim} \psi_1^* p_{\GG}^* T\GG = h_1^* T\GG$:
\begin{align*}
w \in V_q \overline{Q} \mapsto \big( q, dh_1(w) \big),
\end{align*}
noting $dh_1(w) \in T_{h_1(q)} \GG$.  We pull this isomorphism back to $P$ via $s$ to obtain an isomorphism $\alpha_1 : s^* V\overline{Q}|_{X_1} \to s^* h_1^* T\GG|_{X_1} = (h_1 \circ s)^* T\GG|_{X_1}$:
\begin{align*}
(p,w) \mapsto \big( p, dh_1(w) \big).
\end{align*}
Note that $h_1 \circ s$ is a map $P|_{X_1} \to \GG$, and restricting to $X_{01}$ we have a pair of isomorphisms
\begin{align*}
\xymatrixrowsep{2pc}\xymatrixcolsep{1pc}
\xymatrix{
& s^* V\overline{Q}|_{X_{01}} \ar[dl]_{\alpha_0} \ar[dr]^{\alpha_1} & \\
P|_{X_{01}} \times \g & & (h_1 \circ s)^* T\GG. }
\end{align*}
By equivariance of $p_{\GG} \circ \psi_1$, we see that
\begin{align*}
\alpha_1 \circ \alpha_0^{-1} (p, \xi) & = \big( p, dh_1 ( w_\xi(s(p))) \big) = \big( p, u_\xi( h_1 \circ s(p) ) \big) = \big( p, u_\xi ( \sigma(\pi(p)) \cdot g_1(p)) \big).
\end{align*}

Observe that $(h_1 \circ s)^* T\GG$ carries the (right) action
\begin{align*}
(p,w) \cdot g = \big(pg, dR_{g^{-1}} (w) \big),
\end{align*}
and $P|_{X_0} \times \g$ the action
\begin{align*}
(p,\xi) \cdot g = (p \cdot g, \Ad \, g^{-1} \, \xi)
\end{align*}
and $\alpha_0, \alpha_1$ are equivariant with respect to these actions and that on $s^* V\overline{Q}$.

To obtain trivializations for $E = s^* V \overline{Q}/G$, we use the sections $a_i$ of $P$ over $X_i$ to normalize the maps above so that the $G$-component is $e$.  Thus, we take $\beta_0^{-1} : X_0 \times \g \to E|_{X_0}$ as
\begin{align} \label{beta0}
\beta_0^{-1}(y,\xi) & = \big[ a_0(y) : w_\xi \big( b_0(y) \big) \big],
\end{align}
noting that $s \circ a_0 = b_0$; the square brackets indicate the equivalence class modulo the action of $G$.  We wish to define $\beta_1 : E|_{X_1} \to \sigma^* T\GG$ so that
\begin{align*}
\beta_1 \big( [a_1(y) : w] \big) = \big(y, dh_1(w) \big),
\end{align*}
so as to accord with our map $\alpha_1$ above.  Since $p = a_1( \pi(p)) g_1(p)$,
\begin{align*}
[p : w] = [ a_1(\pi(p)) g_1(p) : w] = [a_1(\pi(p)) : dR_{g_1(p)}(w) ],
\end{align*}
so we set
\begin{align} \label{beta1}
\beta_1 \big( [p : w] \big) = \big( \pi(p), dR_{g_1(p)} \circ dh_1(w) \big).
\end{align}
If $w \in V_{s(p)} \overline{Q}$, then $dh_1(w) \in T_{h_1 \circ s(p)} \GG = T_{\sigma(\pi(p)) g_1(p)} \GG$, and so $dR_{g_1(p)} dh_1(w) \in T_{\sigma(\pi(p))} \GG$, so the right side indeed lies in $\sigma^* T\GG$.  It is straightforward to verify that this definition is independent of choice of representative $[p : w]$.

Recalling that $\sigma$ is a map $X_1 \to \GG$, we obtain isomorphisms
\begin{align*}
\xymatrixrowsep{2pc}\xymatrixcolsep{1pc}
\xymatrix{
& E|_{X_{01}} \ar[dl]_{\beta_0} \ar[dr]^{\beta_1} & \\
X_{01} \times \g & & \sigma^* T\GG. }
\end{align*}
Since $X_1$ is not compact $\sigma^* T\GG$ is actually trivial.  However, we will not write it explicitly as a product, but still refer to $\beta_1$ as a trivialization.  The composition of these maps yields
\begin{align} \label{EPQtf}
\begin{split}
\beta_1 \circ \beta_0^{-1}(y,\xi) & = \beta_1 \big ( \big[ a_0(y) : w_\xi \big( b_0(y) \big) \big] \big) = \big( \pi \circ a_0(y), dR_{g_1 \circ a_0(y)} \circ dh_1( w_\xi( b_0(y))) \big) \\
& = \big( y, dR_{g_{10}(y)} u_\xi( h_{10}(y) ) \big) = \big(y, u_{\Ad g_{10}(y) \xi} ( \sigma(y) ) \big).
\end{split}
\end{align}
To avoid overly cumbersome notation later on, we will write this as
\begin{align*} 
\beta_0 \circ \beta_1^{-1} \big( y, u_\xi( \sigma (y) ) \big) = \big( y, \Ad g_{01}(y) \xi \big).
\end{align*}
We may recall that since $\sigma(X_{01}) \subseteq G$, any element of $\sigma^* T\GG$ over $X_{01}$ is of the form $(y, u_\xi( \sigma (y) ) )$ for some $\xi \in \g$.

We will take trivializations $\mu_i : \ad \, P|_{X_i} \to X_i \times \g$ by
\begin{align*}
\mu_i^{-1} (y, \xi) = \big[ a_i(y) : v_\xi \big( a_i(y) \big) \big] = \big[ a_i(y) : \xi \big].
\end{align*}
Then in these trivializations, the map $ds/G : \ad \, P \to E_{PQ}$ looks like
\begin{align} \label{dsGtriv}
\begin{split}
\beta_0 \circ ds/G \circ \mu_0^{-1}(y,\xi) & = (y,\xi) \\
\beta_1 \circ ds/G \circ \mu_1^{-1}(y,\xi) & = \big( y, u_\xi(\sigma(y)) \big).
\end{split}
\end{align}
Note that the first map takes $X_0 \times \g \to X_0 \times \g$ and the second $X_1 \times \g \to \sigma^* T\GG$.

\subsubsection{Infinitesimal Deformations}

We have constructed a map $ds/G : \ad \, P \to E_{PQ}$ of sheaves of sections of vector bundles of the same rank which is an isomorphism except at $x$.  Therefore, we have an exact sequence of sheaves
\begin{align} \label{basseqpb}
0 \to \ad \, P \to E_{PQ} \to E_{PQ}/ \ad \, P \to 0,
\end{align}
and $E_{PQ}/\ad \, P$ is a torsion sheaf.

Since $E_{PQ}/ \ad \, P$ is supported at $x$, representatives (in $E_{PQ}$) of its sections need only be defined in a neighbourhood of $x$, for example, on $X_1$.  But here, one will observe that the expression for $ds/G$ in the $X_1$-trivialization (\ref{dsGtriv}) is precisely that of (\ref{inftmlrightaction}).  This makes clear the following.

\begin{lem} \label{sectioncomparison}
There is an isomorphism
\begin{align*}
H^0(X, E_{PQ}/ \ad \, P) \cong \Gamma(X_1, \sigma^*T\GG) / \Gamma(X_1, \OO^\g_{X_1}).
\end{align*}
If $P$ corresponds to $\varsigma \in \Gr_Q^{\lambda^\vee}$, then this allows us to make the identification
\begin{align*}
T_\varsigma \Gr_Q^{\lambda^\vee} = H^0(X, E_{PQ}/ \ad \, P).
\end{align*}
\end{lem}

\subsubsection{Comparison with the Dual Construction}

We now consider the dual construction to $E_{PQ}$ given in Section \ref{dualconstruction}; we obtain a bundle $E_{QP} := t^* V\overline{P} / G$ which will have transition function (cf.\ (\ref{EPQtf}))
\begin{align} \label{EQPtf}
(y, \xi) \mapsto \big( y, u_{\Ad \, h_{10}(y) \xi} ( \sigma(y)^{-1}) \big).
\end{align}
We claim that $E_{QP} \cong E_{PQ}$.  For this, we use the following.

\begin{lem}
There is an isomorphism $d\tilde{\iota} : \sigma^* T\GG \xrightarrow{\sim} \tau^* T\GG$ under which
\begin{align*}
u_{\Ad \, g_{10}(y) \xi}( \sigma(y) ) \mapsto u_{\Ad \, h_{10} (y) \xi}( \tau(y) ).
\end{align*}
\end{lem}

\begin{proof}
To construct this, we use the involution $\iota : \GG \to \GG$, extending the inversion map on $G$ (Proposition \ref{Gbarinvextn}).  This induces an involution on the tangent bundle $d\iota : T\GG \to T\GG$.  We define $d \tilde{\iota} : \sigma^* T\GG \to \tau^* T\GG$ by
\begin{align*}
(y,v) \mapsto \big( y, -d\iota(v) \big).
\end{align*}
Note that $v \in T_{\sigma(y)}\GG$ and so $-d\iota(v) \in T_{\iota(\sigma(y))}\GG = T_{\tau(y)}\GG$, so this makes sense.  This map is clearly an involution.  Verification of the mapping property is unproblematic.
\end{proof}

\begin{prop}
There is an isomorphism
\begin{align*}
M : E_{PQ} \xrightarrow{\sim} E_{QP}.
\end{align*}
\end{prop}

\begin{proof}
From the above calculations, using the trivializations $\beta_i$ (\ref{beta0}, \ref{beta1}) for $E_{PQ}$ and the analogous trivializations $\gamma_i$ for $E_{QP}$, we define a map $M : E_{PQ} \to E_{QP}$ by
\begin{align} \label{Mtriv}
u \mapsto \begin{cases}
\gamma_0^{-1} \circ \beta_0(u) & u \in E_{PQ}|_{X_0} \\
\gamma_1^{-1} \circ d \tilde{\iota} \circ \beta_1(u) & u \in E_{PQ}|_{X_1}.
\end{cases}
\end{align}
Using (\ref{EPQtf}) and (\ref{EQPtf}) shows that $M$ is well-defined, and it is straightforward to construct an inverse.  This proves that $E_{QP} \cong E_{PQ}$.
\end{proof}

Since the isomorphism $\tilde{\kappa} : \g^* \to \g$ induced by the Killing form $\kappa$ is $\Ad$-equivariant, it induces isomorphisms
\begin{align*}
(\ad \, P)^* & \to \ad \, P, & (\ad \, Q)^* & \to \ad \, Q,
\end{align*}
which we will also denote by $\tilde{\kappa}$; explicitly, if $p \in P, q \in Q, f \in \g^*$,
\begin{align*}
[p : f] & \mapsto [p : \tilde{\kappa}(f)],  & [q : f] & \mapsto [q : \tilde{\kappa}(f)].
\end{align*}
We may now consider the compositions
\begin{align} \label{viaP}
\xymatrix{ E_{PQ}^* \ar[r]^{(ds/G)^t} & (\ad \, P)^* \ar[r]^{ \tilde{\kappa} } & \ad \, P \ar[r]^{ds/G} & E_{PQ} }
\end{align}
and
\begin{align} \label{viaQ}
\xymatrix{ E_{PQ}^* \ar[r]^{(M^{-1})^t}_\sim & E_{QP}^* \ar[r]^{(dt/G)^t} & (\ad \, Q)^* \ar[r]^{ \tilde{\kappa} } & \ad \, Q \ar[r]^{dt/G} & E_{QP} \ar[r]^{M^{-1}} & E_{PQ}, }
\end{align}
and thus get a square
\begin{align} \label{pqqpsquare}
\vcenter{\vbox{
\xymatrix{ E_{PQ}^* \ar[r] \ar[d] & \ad \, Q \ar[d] \\
\ad \, P \ar[r] & E_{PQ}. }
}}
\end{align}

Verifying that it is indeed commutative involves tracing through the maps in the trivializations for the various bundles given above.  In the following, we will denote by a tilde the trivializations for a dual bundle induced by the trivializations on the original bundle.  We wish to check that (\ref{viaP}) and (\ref{viaQ}) yield the same map, which amounts to verifying
\begin{align*}
\beta_i \circ ds/G \circ \tilde{\kappa} \circ (ds/G)^t \circ \tilde{\beta}_i^{-1} = \beta_i \circ M^{-1} \circ dt/G \circ \tilde{\kappa} \circ (dt/G)^t \circ (M^{-1})^t \circ \tilde{\beta}_i^{-1}
\end{align*}
for $i = 0,1$.  It is easy to see from (\ref{dsGtriv}) and (\ref{Mtriv}) that in the $X_0$ trivialization, these maps are all given by the identity on the $\g$ factor, so we need only check the equality in the $X_1$ trivialization.  Using (\ref{dsGtriv}), we may compute the left side acting on $(y, f) \in \sigma^*T\GG$ as
\begin{align} \label{LHScomm}
\beta_1 \circ ds/G \circ \tilde{\kappa} \circ (ds/G)^t \circ \tilde{\beta}_1^{-1} (y,f) = \big( y, dR_{\sigma(y)} \circ \tilde{\kappa} \circ dR_{\sigma(y)}^t(f) \big).
\end{align}
In a similar manner, if we use trivializations of $\ad \, Q$ analogous to the $\mu_i$ for $\ad \, P$, as well as (\ref{Mtriv}), the right side acting on $(y,f)$ becomes
\begin{align} \label{RHScomm}
\big( y, d\iota \circ dR_{\tau(y)} \circ \tilde{\kappa} \circ dR_{\tau(y)}^t \circ d\iota^t(f) \big).
\end{align}
Therefore, it suffices to show that (\ref{LHScomm}) and (\ref{RHScomm}) are equal.  But $\iota \circ R_{\tau(y)} = L_{\sigma(y)} \circ \iota$ and since $\iota$ restricts to the inversion map on $G$, $d\iota : T_e \GG = \g \to T_e \GG = \g$ is $-\1$, so this equality can be readily checked using Corollary \ref{lrtranspose}.  This completes the verification that (\ref{pqqpsquare}) commutes.

\subsection{Calculations for Cocharacters}

We now give a concrete description of the space of infinitesimal deformations $H^0(X, E_{PQ}/ \ad \, P)$ which, by Lemma \ref{sectioncomparison}, amounts to describing the quotient space
\begin{align*}
\Gamma(X_1, \sigma^* T\GG)/ \Gamma(X_1, \OO^\g_{X_1}).
\end{align*}
As in Section \ref{openaffine}, we can use $1/\alpha_1, \ldots, 1/\alpha_l$ to coordinatize $T$.  Since $Y(T) = \Lambda^\vee$, $1/\lambda_1^\vee, \ldots, 1/\lambda_l^\vee$ give a dual basis of cocharacters of $T$.  Thus, if (writing additively now) $\lambda^\vee = - \sum_{i=1}^l r_i \lambda_i^\vee$ is an arbitrary cocharacter, then we may write
\begin{align} \label{cochar}
\lambda^\vee(z) = ( z^{r_1}, \ldots, z^{r_l}).
\end{align}

As noted in Section \ref{HMPBmodfn}, by choosing our trivializations $\psi_1, \varphi_1$ and a coordinate $z$ on $X_1$ centred at $x$ appropriately, we may assume $\sigma$ is a cocharacter. We will thus assume that such choices have been made and that $\sigma$ is of the form (\ref{cochar}) above.  Using the notation of Section \ref{openaffine}, we are setting $z_i = z^{r_i}$:
\begin{align*}
\sigma(z) = F(z^{r_1}, \ldots, z^{r_l}) = \left[ \varphi_0 + \sum_{i=1}^l z^{r_i} \varphi_i + \sum_{k > l} z^{\sum_{i=1}^l r_i n_{ik} } \varphi_k \right].
\end{align*}
Recall that the map $\OO^\g_{X_1} \to \sigma^* T\GG$ is given by (\ref{inftmlrightaction}). Therefore, the image of $\Gamma(X_1, \OO^\g_{X_1})$ consists of all those sections of the form
\begin{align*}
dL(\xi) = dL_{\sigma(z)}(\xi)
\end{align*}
with $\xi \in \g$.  Now, the sections
\begin{align*}
& dA_{(e,e,\sigma(z))}(x_\alpha), & & dA_{(e,e,\sigma(z))}(y_\alpha), & & dA_{(e,e,\sigma(z))}( e_j)
\end{align*}
give a trivialization of $\sigma^* T\GG$, so we want to consider these modulo the image.  But from Lemma \ref{infaction2}, we see that
\begin{align} \label{Tvaluedsigma}
dL(x_\alpha) & = -dA(x_\alpha), & dL(y_\alpha) & = \prod_{i=1}^l z^{ r_i \langle \lambda_i^\vee, \alpha \rangle} & dL(h_i) & = - \sum_{j=1}^l a_{ij} z^{r_j} \, dA(e_j),
\end{align}
recalling that $(a_{ij})$ is the Cartan matrix.  By inverting it, we can see that the sections
\begin{align*}
z^{r_i} dA(e_i), \quad 1 \leq i \leq l
\end{align*}
lie in the image.  Therefore, we can take the sections
\begin{align*} 
\begin{cases} z^j dA(y_\alpha), & 0 \leq j \leq \sum r_i \langle \lambda_i^\vee, \alpha \rangle - 1, \quad \alpha \in \Phi^+, \\ z^k dA(e_i), & 0 \leq k \leq r_i-1, \quad 1 \leq i \leq l,
\end{cases}
\end{align*}
as representative sections for the quotient.  We get analogous expressions if instead we consider the cocharacter $\nu \cdot \lambda^\vee$ with $\nu \in W$.

To simplify the counting of the dimensions, we will take $\sigma = -\lambda_i^\vee$ to be of the form $\sigma(z) = F(1, \ldots, z, \ldots, 1)$, with the $z$ in the $i$th position, so that $r_j = 0$ for $j \neq i$ and $r_i = 1$.  Thus $\sigma$ gives a simple modification.  In this case, the quotient will have representative sections
\begin{align*}
\begin{cases}
z^j dA(y_\alpha), & 0 \leq j \leq \langle \lambda^\vee_i, \alpha \rangle - 1, \quad \alpha \in \Phi^+, \\
dA(e_i). &  \end{cases}
\end{align*}
Therefore, the dimension of the deformation space is
\begin{align*}
\sum_{\alpha \in \Phi^+} \langle \lambda^\vee_i, \alpha \rangle + 1 = 2 \langle \lambda^\vee_i, \rho \rangle + 1,
\end{align*}
which is precisely what we obtained for the dimension of the space of Hecke modifications of type $\lambda^\vee_i$ (\ref{dimYQlambda}).

\subsection{Multiple Modifications} \label{multiplePB}

Consider now a family of bundles given by introducing Hecke modifications of types $\mu_1^\vee, \ldots, \mu_m^\vee$, say with a modification of type $\mu_i^\vee$ at the distinct points $x^i_1, \ldots, x^i_{k_i}$.  Then from Proposition \ref{symmetricproductbundle}, the appropriate parameter space is
\begin{align*}
\left( \Gr_{Q,0}^{\mu_1^\vee, \ldots, \mu_m^\vee} \right)^{(k_1, \ldots, k_m)}.
\end{align*}
%As mentioned in \S\ref{symmetricproducts}, this carries a natural projection to $X^{(k)}$, where $k := \sum_{i=1}^m k_i$, with image the set of reduced divisors in $X^{(k)}$.  We obtain the following.

\begin{prop} \label{multspacePB}
Suppose $P$ is obtained from $Q$ by a sequence of Hecke modifications at distinct points $x^1_1, \ldots, x^1_{k_1}, x^2_1, \ldots, x^2_{k_2}, \ldots, x^m_1, \ldots, x^m_{k_m}$ with a modification of type $\mu^\vee_i$ at $x^i_j, 1 \leq j \leq k_i$, and suppose this corresponds to
\begin{align*}
\varsigma \in \left( \Gr_{Q,0}^{\mu_1^\vee, \ldots, \mu_m^\vee} \right)^{(k_1, \ldots, k_m)}.
\end{align*}
Let $s_\varsigma : P \to \overline{Q}$ be the corresponding map of bundles.  Then one may make the identification
\begin{align*}
T_\varsigma \left( \Gr_{Q,0}^{\mu_1^\vee, \ldots, \mu_m^\vee} \right)^{(k_1, \ldots, k_m)} & = H^0( X, E_{PQ}/ \ad \, P).
\end{align*}

The Kodaira--Spencer map for the family of bundles provided by Proposition \ref{symmetricproductbundle} is given by the connecting homomorphism for the short exact sequence (\ref{basseqpb}):
\begin{align*}
H^0( X, E_{PQ}/ \ad \, P) \to H^1( X, \ad \, P).
\end{align*}
\end{prop}

To see why the last statement holds, observe that a section of $H^0( X, E_{PQ}/ \ad \, P)$ can be represented by a section of $E$ over an open set containing the support of the Hecke modifications, e.g.\ the union of discs $X_1$ that was used above.  In terms of the trivializations described in Section \ref{trivializationsEPQ}, the inclusion $\ad \, P \to E_{PQ}$ is precisely the map in (\ref{localdefmnmap}), and $H^0( X, E_{PQ}/ \ad \, P)$ can be identified with the group in (\ref{localdefmn}).  Recall that in terms of \v{C}ech representatives, the connecting homomorphism is obtained by applying the \v{C}ech differential to such a lift, so in our case, we are restricting this lift to the set $X_{01}$, a disjoint union of punctured discs on which $\ad \, P$ is isomorphic to $E_{PQ}$, so we may think of this as a section of $\ad \, P$ over $X_{01}$, and hence we get a $1$-cocycle with values in $\ad \, P$.  This is at the same time the image of the connecting homomorphism and precisely the deformation of $P$ induced by the section of $E_{PQ}/ \ad \, P$ we started with.

\section{Parametrization of the Moduli Space} \label{Parametrization}

\subsection{Outline}

The moduli space of stable bundles of topological type $\varepsilon \in \pi_1(G)$ will be written $\NN_{G, \varepsilon}$, omitting subscripts when there is no ambiguity; the moduli space of (equivalence classes of) semi-stable bundles will be written $\NN_{G, \varepsilon}^{ss}$.  We now outline our general procedure for constructing parametrizations of the moduli spaces of stable bundles to give a framework for the detailed calculations for specific groups in the next section.  We will start with a fixed bundle $Q$ and consider a family of bundles as in Section \ref{multiplePB} parametrized by
\begin{align*} 
\left( \Gr_{Q,0}^{\mu_1^\vee, \ldots, \mu_m^\vee} \right)^{(k_1, \ldots, k_m)}.
\end{align*}
Breaking with previous notation, we will write $P = P^\varsigma$ for the bundle corresponding to the sequence $\varsigma \in ( \Gr_{Q,0}^{\mu_1^\vee, \ldots, \mu_m^\vee} )^{(k_1, \ldots, k_m)}$ of Hecke modifications of $Q$ at distinct points.  If specification is required, we will write $E_\varsigma$ for $E_{P^\varsigma Q}$.  We will write $N$ for the dimension of this parameter space:
\begin{align*}
N := \dim \left( \Gr_{Q,0}^{\mu_1^\vee, \ldots, \mu_m^\vee} \right)^{(k_1, \ldots, k_m)} = \sum_{i=1}^m k_i \left( 2 \langle \mu_i^\vee, \rho \rangle + 1 \right) = h^0(X, E_\varsigma / \ad \, P),
\end{align*}
the last by Proposition \ref{multspacePB}.  We should observe that by repeated application of Proposition \ref{HMtoptype}, the bundles in this family will have topological type
\begin{align} \label{Ptoptype}
\varepsilon := \varepsilon(Q) + \sum_{i=1}^m k_i [\mu_i^\vee] \in \pi_1(G).
\end{align}

Recall that to $\varsigma \in ( \Gr_{Q,0}^{\mu_1^\vee, \ldots, \mu_m^\vee} )^{(k_1, \ldots, k_m)}$, there corresponds an equivariant map $s_\varsigma : P^\varsigma \to \overline{Q}$.  Since we are only interested in the isomorphism classes of $P^\varsigma$ that arise, rather than the maps $s$ themselves, we will want to quotient out by the action of $\Aut \, \overline{Q} = \Aut \, Q$ (Lemma \ref{AutQbar}).  Therefore we wish to form families of
\begin{align*}
N = \dim \NN_G + \dim \Aut \, Q = \dim \, G (g-1) + h^0(X, \ad \, Q)
\end{align*}
parameters.

We will fix a $P$ in such a family.  From our basic deformation sequence (\ref{basseqpb}), we obtain an exact sequence
\begin{align*}
H^0(X, E) \to H^0(X, E/ \ad \, P) \to H^1(X, \ad \, P) \to H^1(X, E) \to 0
\end{align*}
and the connecting homomorphism is the Kodaira--Spencer map.  In order for us to obtain a parametrization, the family must necessarily be locally complete and so this map must be surjective.  This is equivalent to
\begin{align*}
H^1(X, E) = 0.
\end{align*}
In such a case,
\begin{align*}
h^0(X, E) & = \chi(E) = \chi( E/ \ad \, P) + \chi( \ad \, P) = N + \dim \, G (1-g) = h^0(X, \ad \, Q).
\end{align*}
Since $\ad \, Q \to E$ is an injective map of sheaves, it follows that
\begin{align*} 
H^0(X, \ad \, Q) \to H^0(X, E)
\end{align*}
is an isomorphism (it is injective and we have just seen that both spaces are of the same dimension).  Extending (\ref{pqqpsquare}) into a morphism of short exact sequences, we see that there is a commutative square
\begin{align*}
\xymatrix{ H^0(X, \ad \, Q) \ar[r] \ar[d] & H^0(X, \ad \, Q/ E^*) \ar[d] \\
H^0(X, E) \ar[r] & H^0(X, E/ \ad \, P), }
\end{align*}
it then follows that image of the infinitesimal automorphisms $H^0(X, \ad \, Q)$ of $Q$ in the infinitesimal deformation space $H^0(X, E/ \ad \, P)$ is precisely the kernel of the Kodaira--Spencer map, and therefore the Kodaira--Spencer map
\begin{align*}
H^0(X, E/\ad \, P)/ H^0(X, \ad \, Q) \to H^1(X, \ad \, P)
\end{align*}
for our effective parameter space (i.e., for the quotient by $\Aut \, Q$) is an isomorphism.

By a semicontinuity argument, if $H^1(X, E_\varsigma) = 0$, then $H^1(X, E_\vartheta) = 0$ for all $\vartheta$ in an open neighbourhood of $\varsigma$.  We then get a submersion from an open set in $( \Gr_{Q,0}^{\mu_1^\vee, \ldots, \mu_m^\vee} )^{(k_1, \ldots, k_m)}$ onto one in $\NN_G^{ss}$ and hence the family of bundles obtained must contain stable ones.  Indeed, as this is the case, the modifications yielding stable bundles form a Zariski open set in the parameter space \cite[Proposition 4.1]{Ramanathanthesis} and hence we get an isomorphism of an open set
\begin{align*}
\Omega \subseteq \left( \Gr_{Q,0}^{\mu_1^\vee, \ldots, \mu_m^\vee} \right)^{(k_1, \ldots, k_m)}
\end{align*}
modulo the action of $\Aut \, Q$ onto an open set in $\NN_G$.  To recap, we make the following statement.

\begin{prop} \label{paramznprop}
Suppose we construct a family of bundles via a sequence of Hecke modifications of $Q$ such that the dimension of the parameter space of the family is
\begin{align*}
\dim \, G (g-1) + \dim \Aut \, Q.
\end{align*}
Then if some $P^\varsigma$ in this family is such that $H^1(X, E_{PQ}) = 0$, then in a neighbourhood of $\varsigma$, we obtain a parametrization of $\NN_G$.
\end{prop}

\subsection{Bundles Reducible to a Torus}

Here we will outline a method of constructing families satisfying the hypotheses of Proposition \ref{paramznprop}.  We do this by starting with $Q$ being the trivial bundle and introducing modifications taking values in $T$; this is done to facilitate the calculation of $H^1(X, E)$.  Therefore, from Proposition \ref{paramznprop}, we will want a family of
\begin{align*}
N = \dim \, G \cdot g
\end{align*}
parameters.  The bundle $P$ that we obtain will be reducible to $T$, and hence will not be stable, but in an open neighbourhood of this $P$ in the moduli space, there will still be a Zariski open set consisting of stable bundles.  This is why we need to appeal to the semicontinuity argument mentioned above, to ensure that $H^1(X, E)$ also vanishes for these bundles.

From the root space decomposition of $\g$, it follows that $\ad \, Q$ has a decomposition
\begin{align*}
\ad \, Q = \OO^\g = \OO^\t \oplus \bigoplus_{\alpha \in \Phi} \OO^{\g_\alpha}.
\end{align*}
Let $x \in X$ and $X_0, X_1$ be as before, and consider a modification at $x$ such that with respect to the trivialization of $Q$ and a suitable choice of coordinate $z$ at $x$ it can be represented by the cocharacter $\sigma = \lambda^\vee$ of (\ref{cochar}).  Then this is the transition function $g_{01}$ for a $T$-bundle $R$ and we may write $P = R \times_T G$ as the induced bundle.  One has a decomposition
\begin{align*}
\ad \, P = \OO^\t \oplus \bigoplus_{\alpha \in \Phi} (\ad \, P)_\alpha.
\end{align*}
where $(\ad \, P)_\alpha := R \times_{\Ad \, T} \g_\alpha$.

We observe
\begin{align*}
\alpha \big( \sigma(z) \big) = z^{\langle \lambda^\vee, \alpha \rangle }
\end{align*}
and hence these are the transition functions for $(\ad \, P)_\alpha$ so it follows that
\begin{align*}
(\ad \, P)_\alpha = \OO( - \langle \lambda^\vee, \alpha \rangle x ).
\end{align*}

From the relations (\ref{Tvaluedsigma}), if $\alpha \in \Phi^+$ (respectively, $\alpha \in \Phi^-$) and we let $L_\alpha \subseteq E$ be the line bundle spanned by $x_\alpha$ (respectively, $y_\alpha$) in the $X_0$-trivialization and $dA(x_\alpha)$ (respectively, $dA( y_\alpha)$) in the $X_1$-trivialization, then we see that $E$ is of the form
\begin{align*}
E = V \oplus \bigoplus_{\alpha \in \Phi} L_\alpha,
\end{align*}
where $V$ is a rank $l$ vector bundle.  We will often refer to the $L_\alpha$ as the root bundles of $E$.  Since our goal is to arrange for $H^1(X,E) = 0$, it is clearly enough to show that $H^1(X, L_\alpha) = 0$ for each $\alpha \in \Phi$ and that $H^1(X, V) = 0$.

Under $ds/G : \ad \, P \to E$, we have
\begin{align*}
\OO^\t & \to V, & (\ad \, P)_\alpha & \to L_\alpha, \quad \alpha \in \Phi.
\end{align*}
Indeed, (\ref{Tvaluedsigma}) tells us what the maps look like in the $X_1$-trivialization:
\begin{align} \label{Hmmap}
h_i = \alpha_i^\vee & \mapsto - \sum_{j=1}^l a_{ij} z^{r_j} dA(e_j), & x_\alpha & \mapsto -dA(x_\alpha), & y_\alpha & \mapsto z^{- \langle \lambda^\vee, \alpha \rangle} dA(y_\alpha).
\end{align}
We see that if $\alpha \in \Phi^+$, then
\begin{align} \label{LE}
L_\alpha & \cong (\ad \, P)_\alpha = \OO(- \langle \lambda^\vee, \alpha \rangle x), & L_{-\alpha} & \cong (\ad \, P)_{-\alpha} ( \langle \lambda^\vee, \alpha \rangle x) = \OO.
\end{align}

If instead of $\lambda^\vee$, we use the modification corresponding to $\nu \lambda^\vee$, i.e.\ twisted by an element of the Weyl group, then $L_\alpha$ will be spanned by $x_\alpha$ in the $X_0$-trivialization and by $dA^\nu( \Ad \, w^{-1} x_\alpha)$ in the $X_1$-trivialization (cf.\ Lemma \ref{infactiontwist}).  Replacing $x_\alpha$ by $y_\alpha$ as appropriate for the negative root spaces, we see that for $\alpha \in \nu \Phi^+$,
\begin{align} \label{LEtwist}
L_\alpha & \cong (\ad \, P)_\alpha = \OO (- \langle \lambda^\vee, \nu^{-1} \alpha \rangle x), &  L_{-\alpha} & \cong (\ad \, P)_{-\alpha} ( \langle \lambda^\vee, \nu^{-1} \alpha \rangle x) = \OO.
\end{align}

Now, assume that we introduce modifications at the points $x_1, \ldots, x_M$, with the modification at $x_j$ being the cocharacter $-\nu_j \lambda_{i(j)}, 1 \leq j \leq M$.  Then it follows from (\ref{LE}) and (\ref{LEtwist}) that if we define the effective divisors $D_\alpha$ by
\begin{align*}
D_\alpha := \sum_{j=1}^M n^\alpha_j x_j,
\end{align*}
where
\begin{align*}
n^\alpha_j = \max \{ 0, \langle \lambda_{i(j)}, \nu_j^{-1} \alpha \rangle \} = \begin{cases} 0 & \text{if } \alpha \in \nu_j \Phi^- \\ \langle \lambda_{i(j)}, \nu_j^{-1} \alpha \rangle & \text{if } \alpha \in \nu_j \Phi^+ \end{cases}
\end{align*}
then
\begin{align*}
L_\alpha = \OO( D_\alpha).
\end{align*}

To arrange for $H^1(X, L_\alpha) = 0$, we will need $\deg D_\alpha \geq g$; for a generic choice of $D_\alpha$, in a sense explained more fully below, equality will be sufficient.  In the next section, what we will attempt to do is to start with a number of simple Hecke modifications so that we get a parameter space of dimension $\dim \, G \cdot g$, and then apply various elements of the Weyl group so that the root parameters are distributed among the root spaces in a manner just described to ensure that $H^1(X, L_\alpha) = 0$ for all $\alpha \in \Phi$.

We also need to arrange for $H^1(X, V) = 0$.  As mentioned above, $V$ is obtained by upper Hecke modification of the trivial bundle $\OO^\t$; indeed, since a simple Hecke modification introduces a single toral parameters, it is obtained via a sequence of simple upper Hecke modifications, so that there is an exact sequence
\begin{align} \label{torex}
0 \to \OO^\t \to V \to V / \OO^\t \to 0,
\end{align}
with the quotient a sum of skyscraper sheaves.  If we use the cocharacter $-\lambda_i^\vee$ at the point $x \in X$, then from (\ref{Hmmap}), we can see that the kernel of the fibre at $x$ is spanned by
\begin{align} \label{Hmker}
\xi_i := \sum_{k=1}^l a^{ik} \alpha_k^\vee,
\end{align}
considered as an element of $\t = (\OO^\t)_x$, where $(a^{ik})$ is the inverse of the Cartan matrix.  If instead of $-\lambda_i^\vee$, we use the modification corresponding to $-\nu \lambda_i^\vee$, then the kernel is spanned by $\nu \xi_i$.

\subsection{Choosing Points Generically}

We now make precise what is meant by a ``generic'' choice of divisor and use this to give a sufficient condition for the vanishing of $H^1(X, V)$ in a situation as in (\ref{torex}).

\begin{prop} \label{choosediv}
Let $m \geq 1$.  Then there exists $\tilde{D} \in X^{(m)}$ such that if $(g_1, \ldots, g_k)$ is a partition of $g$ and $x_{i_1}, \ldots, x_{i_k} \in \supp \, \tilde{D}$ and we set
\begin{align*}
D = \sum_{j=1}^k g_j x_{i_j},
\end{align*}
then
\begin{align} \label{H1iszero}
H^1(X, \OO(D)) = 0.
\end{align}
In fact, the set of such $\tilde{D}$ is dense in $X^{(m)}$.  In particular, if $m \geq g$ and $\tilde{D}$ is a reduced divisor of degree $m$ (i.e., a sum of $m$ distinct points), then for any $0 \leq D \leq \tilde{D}$ with $\deg D = g$, (\ref{H1iszero}) holds.
\end{prop}

\begin{proof}
First, consider the partition $(g_1 = g)$ and let $x \in X$, so that $D = g \cdot x$.  Then
\begin{align*}
h^1(X, \OO(D)) = h^0(X, \OO(g \cdot x)) - 1
\end{align*}
and this vanishes except when $x$ is a Weierstrass point, of which there are finitely many.

Now, fix an arbitrary partition $(g_1, \ldots, g_k)$, consider the map $X^k \to X^{(g)}$ given by
\begin{align*}
(x_1, \ldots, x_k) \mapsto D := \sum_{j=1}^k g_j x_j,
\end{align*}
and consider the family of line bundles $\OO(D)$ as $D$ ranges over the image of this map.  Taking $x_1 = \cdots = x_k = x$, and this to be a non-Weierstrass point, we see that there exists a line bundle in this family with vanishing first cohomology.  By the semi-continuity theorem, there is a Zariski open, non-empty and hence dense, subset of $X^k$ for which we have $H^1(X, \OO(D)) = 0$.

Let $m \geq 1$ be given.  If $m < k$, then choosing the $x_{i_j}$ in the support of a degree $m$ divisor and forming $D$ as above, there will necessarily be repeated points, so we consider the maps $X^m \to X^k$
\begin{align*}
(x_1, \ldots, x_m) \mapsto (x_1, \ldots, x_m, x_{r_1}, \ldots, x_{r_{k-m}}),
\end{align*}
where the $r_j \in \{ 1, \ldots, m \}$.  The set of points in $X^m$ for which we have vanishing will be the union of the preimages of the open set in $X^k$ for which we have vanishing, so is open and dense in $X^m$ and hence in $X^{(m)}$.

Assume $m \geq k$.  Then choosing a subset of $k$ points amounts to taking a projection map $X^m \to X^k$ onto a set of $k$ factors.  This time, the sets for which we have vanishing for all such choices is the intersection of the preimages of the open set in $X^k$ constructed above.  But since there are finitely many such projections, this is a finite intersection and again we have an open dense subset in $X^m$ and hence $X^{(m)}$ for which we have vanishing.  Thus, for any $m$ and for any partition, we can find a dense open set of $\tilde{D} \in X^{(m)}$ for which the conclusion holds.

Finally, for a fixed $m \geq 1$, since there are only finitely many partitions of $g$, we can take the intersection of the open sets just constructed for each partition and the statement follows.
\end{proof}

The next statement follows easily from the above.

\begin{cor}
If $m \geq g$, then for a generic choice of $m$ points in $X$, no $g$ of them will form part of a canonical divisor.
\end{cor}

\begin{prop} \label{H1Vvanish}
Let $v_1, \ldots, v_l \in \C^l$ be a basis of $\C^l$.  Choose points $x_{ij} \in X, 1 \leq i \leq l, 1 \leq j \leq g$ such that no $g$ of them form part of a canonical divisor (by the preceding Corollary, a generic choice of points satisfies this condition), and suppose $V$ is the vector bundle obtained from the trivial bundle by the upper Hecke modification determined by (the subspace spanned by) $v_i$ at $x_{ij}$ for $1 \leq i \leq g$.  Then
\begin{align*}
H^1(X, V) = 0.
\end{align*}
\end{prop}

\begin{proof}
In (\ref{torex}), if we replace $\OO^\t$ by $\OO^{\oplus l}$, we obtain an exact sequence
\begin{align*}
H^0(X, V/ \OO^{\oplus l}) \cong \C^{gl} \to H^1(X, \OO^{\oplus l}) = \big( H^0( X, K)^* \big)^{\oplus l} \to H^1(X, V) \to 0.
\end{align*}
Then $H^1(X, V) = 0$ if and only if the connecting homomorphism is surjective, and this holds if and only if the transpose map
\begin{align} \label{inj?}
H^0(X, K)^{\oplus l} \to H^0(X, V/\OO^{\oplus l})^*
\end{align}
is injective.  This is the statement we prove.

Recall that sections of $H^0(X, V/ \OO^{\oplus l})$ are given by collections of $l$-tuples
\begin{align*}
f_{ij} = {}^t(f_{ij}^1, \ldots, f_{ij}^l)
\end{align*}
of functions meromorphic in a neighbourhood of $x_{ij}$ for each $1 \leq i \leq l, 1 \leq j \leq g$, having at most a simple pole along $v_i$.  The map (\ref{inj?}) is given by
\begin{align*}
{}^t \omega = (\omega_1, \ldots, \omega_l) \mapsto \left( (f_{ij}) \mapsto \sum_{i=1}^{l} \sum_{j=1}^g \Res_{x_{ij}} {}^t f_{ij} \omega \right),
\end{align*}
where
\begin{align*}
{}^t f_{ij} \omega = \sum_{k=1}^l f^k_{ij} \omega_k.
\end{align*}

Consider the following basis of $H^0(X, V/ \OO^{\oplus l})$:  let $s_{ij} = ( f(i,j)_{rs})$ where
\begin{align*}
f(i,j)_{rs} = \begin{cases} 0 & \text{if } (i,j) \neq (r,s) \\ z^{-1}_{ij} v_i & \text{if } (i,j) = (r,s). \end{cases}.
\end{align*}
where, $z_{ij}$ is a coordinate centred at $x_{ij}$.

If $\omega \in H^0(X, K)^{\oplus l}$ lies in the kernel of this map, then
\begin{align*}
\Res_{x_{ij}} z_{ij}^{-1} ({}^t v_i) \omega = 0
\end{align*}
for $1 \leq i \leq l, 1 \leq j \leq g$.  Fixing $i$, we see that ${}^t v_i \omega$ is a fixed linear combination of $\omega_1, \ldots, \omega_l$, and this differential vanishes at the $g$ points $x_{i1}, \ldots, x_{ig}$.  By our choice of points, it follows that ${}^t v_i \omega = 0$.  But since $v_1, \ldots, v_l$ are linearly independent, we can solve for each $\omega_k$ to obtain $\omega_1 = \cdots = \omega_l = 0$.  This proves the claimed injectivity.
\end{proof}

\begin{rmk}
One will observe that for the proof to go through one needs only that each ${}^t v_i \omega$ vanishes and that the $v_i$ are a basis.  If in our hypothesis, we allow more points, then this vanishing still holds, and the conclusion of the Proposition goes through.  Also, introducing more Hecke modifications only increases the degrees of the root bundles, which does not affect the vanishing of the first cohomology once it is achieved.  Thus, if we construct a family of large enough dimension, possibly much larger than $\dim \, G \cdot g$, we can always obtains submersions from a parameter space onto the moduli space.
\end{rmk}

\section{Calculations for Specific Groups}

In this section, we show that the conditions of Proposition \ref{paramznprop} are satisfied in several situations:  namely, for the adjoint forms of the semisimple groups with root systems $A_3, C_l, D_l$ when the genus of $X$ is even.  We use the standard labelling of the simple roots for each of the root systems (as used in Chapter 11 of \cite{Humphreys}, for example).  For each Hecke modification corresponding to the negative of a fundamental coweight $-\lambda_i^\vee$, we find the roots which get twisted as in (\ref{LE}); these are precisely the $\alpha \in \Phi^+$ in which $\alpha_i$ appears when written as a sum of simple roots.  We compute the number of parameters obtained for each of these modification types, namely the coefficient of the $\alpha_i$ which is $\langle \lambda_i^\vee, \alpha \rangle$.  Also, we will explain how to use the Weyl group to distribute the parameters among the root bundles.

\subsection{Calculations for $A_l$}

For the root system $A_l$, the positive roots are all of the form $\sum_{k=i}^{j-1} \alpha_k$ for some $1 \leq i < j \leq l+1$.  It is not hard to see that each $\alpha_i$ appears in $i(l+1-i)$ positive roots and each with multiplicity one.  Therefore
\begin{align*}
2 \langle \lambda_i^\vee, \rho \rangle + 1 = i (l+1-i) + 1.
\end{align*}

Because of the symmetry of the root system under
\begin{align*}
\alpha_i \mapsto \alpha_{l-i},
\end{align*}
we need only consider modification types $-\lambda_i^\vee$ for $1 \leq i \leq \lfloor \frac{1}{2}(l+1) \rfloor$.

There is a concrete realization of this root system as the set of vectors
\begin{align*}
\{ \pm ( e_i - e_j) \, | \, 1 \leq i < j \leq l+1 \},
\end{align*}
where $e_1, \ldots, e_{l+1}$ are the standard basis vectors, in the subspace of $\R^{l+1}$ whose coordinates sum to zero.  We can take
\begin{align*}
\alpha_i := e_i - e_{i+1}, \quad 1 \leq i \leq l
\end{align*}
as a set of simple roots.  The Weyl group is the symmetric group $\mathfrak{S}_{l+1}$ and it acts simply by permuting the indices.  The simple reflections $\nu_i$ correspond to the permutations $(i \ i+1)$.

We have only carried out the analysis for $l = 3$.  For $l = 1$, it is possible to obtain a parametrization when $g = 2$, but this requires a slightly different set-up than what we have described here.  For $l = 2$, one can obtain a family with the desired number of parameters for $g = 3$, but it is unclear that we can arrange for $H^1(X, V) = 0$.  For larger $l > 3$, general formulae to obtain families with the correct number of parameters are not apparent.

\subsubsection{$A_3$}

A modification of type $-\lambda_2^\vee$ twists the root bundles indexed by
\begin{align*}
\alpha_2 = e_2 - e_3, & & \alpha_1 + \alpha_2 = e_1 - e_3, & & \alpha_2 + \alpha_3 = e_2 - e_4, & & \alpha_1 + \alpha_2 + \alpha_3 = e_1 - e_4,
\end{align*}
so gives a total of $5$ parameters.  It is easy to check that the permutations $(1 \ 4)$ and $(2 \ 3)$ leave fixed this set of roots.  They generate a subgroup of order 4 in $\mathfrak{S}_4$.  To see which other sets of roots appear in the Weyl orbit, we apply the following coset representatives:
\begin{align*}
(1 \ 3): & & - (e_1 - e_2) & & - (e_1 - e_3) & & e_2 - e_4 & & e_3 - e_4 \\
(2 \ 3): & & - (e_2 - e_3) & & e_1 - e_2 & & e_3 - e_4 & & e_1 - e_4 \\
(1 \ 4): & & e_2 - e_3 & & - (e_3 - e_4) & & -(e_1 - e_2) & & -(e_1 - e_4) \\
(2 \ 4): & & - (e_3 - e_4) & & e_1 - e_3 & & -(e_2 - e_4) & & e_1 - e_2 \\
(1 \ 4) (2 \ 3): & & - (e_2 - e_3) & & - (e_2 - e_4) & & -(e_1 - e_3) & & -(e_1 - e_4).
\end{align*}
One observes that each root appears exactly twice.

Since $\dim \, A_3 = 15$, we would like to construct a family with $N = 15 g$ parameters. We will assume that $g = 2k$ is even.  We take $6k$ modifications of type $-\lambda^\vee_2$, and for each permutation above, we take $k$ of the modifications to be twisted by that permutation.  Then we see that each root bundle receives $2k = g$ parameters, as required.  Since $\pi_1( A_3)$ is cyclic of order 4, and $\lambda_2^\vee$ represents the order 2 element, by (\ref{Ptoptype}), we will get a topologically trivial bundle.

To obtain a parametrization, we must show that $H^1(X, V) = 0$.  From (\ref{Hmker}), for a modification corresponding to $\lambda^\vee_2$ at $x$, the kernel of the map of fibres is spanned by
\begin{align*}
\xi_2 = \tfrac{1}{2} ( \alpha^\vee_1 + 2 \alpha^\vee_2 + \alpha^\vee_3 ).
\end{align*}
Using the correspondence $\nu_i \leftrightarrow (i \ i+1)$, one checks that
\begin{align*}
(1 \ 3 ) \cdot \xi_2 & = - (2 \ 4) \cdot \xi_2 = \tfrac{1}{2} ( - \alpha^\vee_1 + \alpha^\vee_3), & (2 \ 3) \cdot \xi_2 & = - (1 \ 4) \cdot \xi_2 = \tfrac{1}{2} ( \alpha^\vee_1 + \alpha^\vee_3) \\
(1 \ 4) ( 2 \ 3) \cdot \xi_2 & = \xi_2. & &
\end{align*}
Therefore, taking
\begin{align*}
v_1 & = \xi_2, & v_2 & = (1 \ 3) \cdot \xi_2, & v_3 & = (2 \ 3) \cdot \xi_2,
\end{align*}
we get a basis of $\t$ and since each $v_i$ spans the line determining the Hecke modification at $2k = g$ points, the hypotheses of Lemma \ref{H1Vvanish} are satisfied.

\subsection{Calculations for $C_l$}

The root system $C_l$ has positive roots
\begin{align*}
\alpha_i + \cdots + \alpha_{j-1} & & 1 \leq i < j \leq l \\
\alpha_i + \cdots + \alpha_{j-1} + 2\alpha_j + \cdots + 2 \alpha_{l-1} + \alpha_l & & 1 \leq i < j \leq l-1 \\
\alpha_i + \cdots + \alpha_l & & 1 \leq i \leq l \\
2\alpha_i + \cdots + 2\alpha_{l-1} + \alpha_l & & 1 \leq i \leq l-1.
\end{align*}
From this, one finds
\begin{align*}
2 \langle \lambda_i^\vee, \rho \rangle + 1 = \begin{cases} i (2l - i + 1) + 1 & 1 \leq i \leq l-1 \\ \tfrac{1}{2} l (l+1) + 1 & i = l. \end{cases}
\end{align*}

Recall that the root system can be realized as the vectors
\begin{align*}
\{ \pm e_i \pm e_j, \pm 2e_k \in \R^l \, | \, 1 \leq i < j \leq l, 1 \leq k \leq l \}.
\end{align*}
The simple roots are
\begin{align*}
\alpha_i & := e_i - e_{i+1}, \quad 1 \leq i \leq l-1, & \alpha_l & := 2e_l.
\end{align*}

A modification of type $-\lambda_1^\vee$ contains the roots
\begin{align*}
e_1 - e_2, \ldots, e_1 - e_l, \quad e_1 + e_2, \ldots, e_1 + e_l, \quad 2e_1,
\end{align*}
and $\langle \lambda_1^\vee, \alpha \rangle = 1$ for each of these roots $\alpha$, except the last, where we get $2$.  We get a total of $2l + 1$ parameters for each such modification.

The Weyl group of $C_l$ is the semidirect product $(\Z/ 2\Z)^l \rtimes \mathfrak{S}_l$; $(a_1, \ldots, a_l) \in (\Z/ 2\Z)^l, \sigma \in \mathfrak{S}_l$ act as follows
\begin{align*}
(a_1, \ldots, a_l) e_i & = (-1)^{a_i} e_i & \sigma e_i & = e_{\sigma i},
\end{align*}
and we interpret an element $( (a_1, \ldots, a_l), \sigma)$ as the composition of these actions.

We consider the element $\nu = (1, 0, \ldots, 0, (1 \ 2 \cdots l )) \in W$.  It is straightforward to verify that
\begin{align*}
e_i - e_j & = \nu^{i-1} ( e_1 - e_{j-i+1}) = \nu^{l+j-1} ( e_1 + e_{l+i-j+1}) \\
-(e_i - e_j) & = \nu^{l+i-1} ( e_1 - e_{j-i+1}) = \nu^{j-1} ( e_1 + e_{l+i-j+1}) \\
e_i + e_j & = \nu^{i-1} ( e_1 + e_{j-i+1}) = \nu^{j-1} ( e_1 - e_{l+i-j+1}) \\
-(e_i + e_j) & = \nu^{l+i-1} ( e_1 + e_{j-i+1}) = \nu^{l+j-1} ( e_1 - e_{l+i-j+1}) \\
2e_i & = \nu^{i-1} (2e_1) \\
-2e_i & = \nu^{l+i-1} (2e_1).
\end{align*}
Therefore, using the modifications $- \lambda_1^\vee, - \nu \lambda_1^\vee, \ldots, -\nu^{2l-1} \lambda_1^\vee$, we obtain each root with multiplicity two, with the long roots obtained twice with different modifications, while the short roots are obtained with multiplicity two by a single modification.

Recall that $\dim \, C_l = l(2l+1)$, so we want $N = lg(2l+1)$ parameters.  Since a modification of type $-\lambda_1$ yields $2l+1$ parameters, we will take $M = lg$.  Suppose that $g = 2m$ is even.  Then we will put in a modification of type $-\nu^k \lambda_1$ at $m$ points, for $0 \leq k \leq 2l-1$.  As such, if $\alpha$ is a long root, then $L_\alpha$ is of the form $\OO(D)$ where $D$ is a divisor of degree $2m = g$, whose support may be taken to be distinct points.  However, if $\alpha$ is a short root, then $L_\alpha$ is of the form $\OO(2D)$ for some degree $m$ divisor.  In either case, by Proposition \ref{choosediv}, we can choose the points generically so that these have vanishing first cohomology groups.

To show that $H^1(X, V) = 0$, we consider the element
\begin{align*}
\xi_1 = \alpha_1^\vee + \cdots + \alpha_{l-1}^\vee + \tfrac{1}{2} \alpha_l^\vee
\end{align*}
under the action of $\nu$.  It is not hard to see that
\begin{align*}
\nu = \nu_1 \cdots \nu_l,
\end{align*}
where the $\nu_i$ are the simple reflections.  We compute
\begin{align*}
\nu \cdot \alpha_i^\vee & = \begin{cases} \alpha_{i+1}^\vee & i = 1, \ldots, l-2 \\
\alpha_1^\vee + \cdots + \alpha_{l-1}^\vee + 2 \alpha_l^\vee & i = l-1 \\
-( \alpha_1 + \cdots + \alpha_l) & i = l. \end{cases}
\end{align*}
From this, one sees that
\begin{align*}
\nu^i \cdot \xi_i = \tfrac{1}{2} \big( \alpha_i^\vee + 3( \alpha_{i+1}^\vee + \cdots + \alpha_l^\vee) \big),
\end{align*}
for $1 \leq i \leq l-1$.  With a bit of work, one can then show that $\xi_1, \nu \cdot \xi_1, \ldots, \nu^{l-1} \cdot \xi_1$ are linearly independent (in fact, if one writes these vectors in terms of the basis $\alpha_i^\vee$, then the resulting matrix has determinant $(-1)^{l-1} - 2^{-l} \neq 0$).  Observing that $\nu^l = -\1$, with the modifications introduced as above, we can apply Proposition \ref{H1Vvanish}, to see that $H^1(X, V) = 0$.

As to the topological type, we recall that the fundamental group of $C_l$ is of order 2, and since we are using an even number of modifications and starting from the trivial bundle, we again get a parametrization of the moduli space of topologically trivial bundles.

\subsection{Calculations for $D_l$}

The positive roots in $D_l$ are
\begin{align*}
\alpha_i + \cdots + \alpha_{j-1} & & 1 \leq i < j \leq l \\
\alpha_i + \cdots + \alpha_{j-1} + 2\alpha_j + \cdots + 2\alpha_{l-2} + \alpha_{l-1} + \alpha_l & & 1 \leq i < j \leq l-2 \\
\alpha_i + \cdots + \alpha_l & & 1 \leq i \leq l-2, i = l \\
\alpha_i + \cdots + \alpha_{l-2} + \alpha_l & & 1 \leq i \leq l-2 .
\end{align*}
From this, one may readily compute
\begin{align*}
2 \langle \lambda^\vee_i, \rho \rangle + 1 = \begin{cases} i(2l-i-1) + 1 & 1 \leq i \leq l-2 \\ \tfrac{1}{2}l(l - 1) + 1 & i = l-1, l. \end{cases}
\end{align*}

The root system for $D_l$ can be realized as the set of vectors
\begin{align*}
\{ \pm e_i \pm e_j \in \R^l \, | \, 1 \leq i < j \leq l \}.
\end{align*}
A set of simple roots can be given by
\begin{align*}
\alpha_i & = e_i - e_{i+1}, \quad 1 \leq i \leq l-1, & \alpha_l & = e_{l-1} + e_l.
\end{align*}
The Weyl group is $(\Z/2\Z)^{l-1} \ltimes \mathfrak{S}_l$, where $(\Z/2\Z)^{l-1}$ is realized as the subgroup of $(\Z/2\Z)^l$ the sum of whose components is even.  The action is otherwise the same as that for the Weyl group of $C_l$.  Recall that $-\1 \in W$ if and only if $l$ is even; it corresponds to the element $(1, \ldots, 1, e)$.

The positive roots containing $\alpha_1$ are
\begin{align*}
e_1 - e_j, \quad e_1 + e_j & & 2 \leq j \leq l
\end{align*}
and so these are the indices of the root bundles twisted when we introduce a modification of type $-\lambda^\vee_1$; with the toral parameter, we see that a modification of this type yields $2l-1$ parameters.

In the case that $l$ is even, so that $(1, \ldots, 1, e) \in W$, we consider the $2l$ modifications given by
\begin{align*}
& - \big(0, \ldots, 0, (1 \ i) \big) \cdot \lambda^\vee_1, \quad 1 \leq i \leq l, & & - \big(1, \ldots, 1, (1 \ i) \big) \cdot \lambda^\vee_1, \quad 1 \leq i \leq l,
\end{align*}
where for $i = 1$ by $(1 \ i)$ we mean the identity permutation.
Observe that
\begin{align*}
e_i - e_j & = \big( 0, \ldots, 0, (1 \ i) \big) \cdot ( e_1 - e_j) = \big( 1, \ldots, 1, (1 \ j) \big) \cdot ( e_1 - e_i) & & 1 \leq i \leq l \\
-(e_i - e_j) & = \big(1, \ldots, 1, (1 \ i) \big) \cdot ( e_1 - e_j) = \big( 0, \ldots, 0, (1 \ j) \big) \cdot ( e_1 - e_i) & & 1 \leq i \leq l \\
e_i + e_j & = \big( 0, \ldots, 0, (1 \ i) \big) \cdot ( e_1 + e_j) = \big( 0, \ldots, 0, (1 \ j) \big) \cdot ( e_1 + e_i) & & 1 \leq i \leq l \\
-(e_i + e_j) & = \big(1, \ldots, 1, (1 \ i) \big) \cdot ( e_1 + e_j) = \big( 1, \ldots, 1, (1 \ j) \big) \cdot ( e_1 + e_i) & & 1 \leq i \leq l.
\end{align*}
Hence each root appears exactly twice.

In the case that $l$ is odd, if $\nu := (1, \ldots, 1, 0, e) \in W$, we consider the $2l$ modifications given by
\begin{align*}
& - (1 \ i) \cdot \lambda^\vee_1, \quad 1 \leq i \leq l, & & -(1 \ i) \nu \cdot \lambda^\vee_1, \quad 1 \leq i \leq l.
\end{align*}
Here, we can see that
\begin{align*}
e_i - e_j & = (1 \ i) \cdot ( e_1 - e_j) = (1 \ j) \nu \cdot ( e_1 - e_i) & & 1 \leq i < j \leq l \\
-(e_i - e_j) & = (1 \ i) \nu \cdot ( e_1 - e_j) = (1 \ j) \cdot ( e_1 - e_i) & & 1 \leq i < j \leq l \\
e_i + e_j & = (1 \ i) \cdot ( e_1 + e_j) = (1 \ j) \cdot ( e_1 + e_i) & & 1 \leq i < j \leq l \\
-(e_i + e_j) & = (1 \ i) \nu \cdot ( e_1 + e_j) = (1 \ j) \nu \cdot ( e_1 + e_i) & & 1 \leq i < j \leq l .
\end{align*}
Again, every root appears exactly twice.

We will assume that $g = 2k$ is even.  Recall that $\dim D_l = 2l^2 - l = l(2l-1)$, so that we want $N = lg(2l-1) = 2kl(2l-1)$ total parameters.  Since each modification introduces $2l-1$ parameters, we will use $lg = 2kl$ of them.  We have just seen that with $2l$ modifications, we can introduce $2$ parameters to each root bundle, so with $2kl$, we get $2k = g$ in each bundle, as required.

We now show that $H^1(X, V) = 0$.  We consider the orbit of
\begin{align*}
\xi_1 = \sum_{k=1}^l a^{1k} \alpha^\vee_k = \alpha^\vee_1 + \cdots + \alpha^\vee_{l-2} + \tfrac{1}{2} ( \alpha^\vee_{l-1} + \alpha^\vee_l )
\end{align*}
under the elements of the Weyl group used above.  Using this expression, one may check that for $2 \leq i \leq l-3$, we have
\begin{align*}
\nu_i \cdot \xi_1 = \xi_1.
\end{align*}
Since $\nu_i$ corresponds to the permutation $(i \ i+1)$ for $1 \leq i \leq l-1$, the correspondence
\begin{align*}
(1 \ i) \leftrightarrow \nu_{i-1} \nu_{i-2} \cdots \nu_2 \nu_1 \nu_2 \cdots \nu_{i-2} \nu_{i-1},
\end{align*}
for $2 \leq i \leq l$, follows.  So by induction, one obtains
\begin{align*}
(1 \ i) \cdot \xi_1 & = \xi_1 - ( \alpha^\vee_1 + \cdots + \alpha^\vee_{i-1}),
\end{align*}
for $2 \leq i \leq l-2$.  Also,
\begin{align*}
(1 \ l-1 ) \cdot \xi_1 & = \nu_{l-2} \big( \alpha^\vee_{l-2} + \tfrac{1}{2} ( \alpha^\vee_{l-1} + \alpha^\vee_l ) \big) = \tfrac{1}{2} ( \alpha^\vee_{l-1} + \alpha^\vee_l ),  \\
(1 \ l) \cdot \xi_1 & = \nu_{l-1} \cdot \tfrac{1}{2} ( \alpha^\vee_{l-1} + \alpha^\vee_l ) = \tfrac{1}{2} ( -\alpha^\vee_{l-1} + \alpha^\vee_l ).
\end{align*}
From all of this, it is not hard to see that
\begin{align*}\xi_1, (1 \ 2) \cdot \xi_1, \ldots, (1 \ l) \cdot \xi_1\end{align*}
are linearly independent.

Consider the element $\nu = (1, \ldots, 1, 0, e)$ in the case that $l$ is odd.  Then we can see that
\begin{align*}
\nu \cdot \alpha_i & = -\alpha_i, \quad 1 \leq i \leq l-2, & \nu \cdot \alpha_{l-1} & = - \alpha_l, & \nu \cdot \alpha_l & = - \alpha_{l-1}.
\end{align*}
Since the matrix for $\nu$ acting on $\t$ with respect to the basis $\alpha^\vee_i$ is the same as that for its action on $\t^*$ with respect to the basis $\alpha_i$, it follows that
\begin{align*}
\nu \cdot \xi_1 = -\xi_1.
\end{align*}
Therefore, in the case where $l$ is odd, the same argument as above can be used to show that we can arrange for $H^1(X, V) = 0$.

One will recall that $\pi_1(D_l)$ is either $\Z/4\Z$ or $\Z/2\Z \times \Z/2\Z$, depending on the parity of $l$.  However, in either case, $\lambda_1^\vee$ represents an element of order 2, and since we are using an even number of modifications and starting with the trivial bundle, we again get a parametrization of the moduli space of topologically trivial bundles.

\begin{rmk}
There are two obstacles for us in obtaining the parametrizations we seek.  The first is to construct families of the requisite dimension $N = \dim \, G \cdot g$, given the number of parameters yielded by modifications of each type.  As we mentioned, this is a problem for $A_l$ when $l > 3$.  It is also a problem for $B_l$:  a modification of type $-\lambda_i^\vee$ yields $i(2l-i) + 1$ parameters and $\dim \, G = l(2l+1)$.  It is not clear what combination of these types would yield the correct number of parameters.

The second obstacle is to find a way to evenly distribute the parameters among the root bundles using the Weyl action.  In the case of $G_2$, one of the simple Hecke modifications yields $7$ parameters.  Since $\dim \, G_2 = 14$, one can always obtain the desired number of parameters.  However, the problem is that this modification introduces $2$ parameters corresponding to short roots and $4$ corresponding to the long roots, yet there are the same number of short and long roots in $G_2$.  Therefore it is impossible to obtain the required $g$ parameters in the short root spaces with the given number of modifications.
\end{rmk}

\small

\bibliographystyle{amsalpha}
\bibliography{research}

\end{document}